\documentclass[11pt]{article}
\usepackage{amsfonts,amsmath}
\usepackage{amssymb,mathrsfs,latexsym,bm}
\begin{document}
\newcommand{\qed}{\hphantom{.}\hfill $\Box$\medbreak}
\newcommand{\proof}{\noindent{\bf Proof \ }}
\newtheorem{Theorem}{Theorem}[section]
\newtheorem{Lemma}[Theorem]{Lemma}
\newtheorem{Corollary}[Theorem]{Corollary}
\newtheorem{Remark}[Theorem]{Remark}
\newtheorem{Example}[Theorem]{Example}
\newtheorem{Definition}[Theorem]{Definition}
\newtheorem{Construction}[Theorem]{Construction}

\thispagestyle{empty}
%%%%%%%%%%%%%%%%%%%%%%%%%%%%%%%%%%%%%%%%%%%%%%%%%%%%%%%%%%%%%%%%%%%%%%%%
 \renewcommand{\thefootnote}{\fnsymbol{footnote}}

\begin{center}
{\Large\bf Semi-cyclic holey group divisible designs and applications to sampling designs and optical orthogonal codes \footnote{Supported by the Fundamental Research Funds for the Central Universities   grant $2013$JB$Z005$, the Natural Science Foundation of Ningbo grant $2013$A$610102$,  the NSFC  grants $11471032$ and $11201252$, and NSERC grant 239135-2011}}

\vskip12pt

Tao Feng$^1$, Xiaomiao Wang$^2$ and Ruizhong Wei$^3$ \\[2ex] {\footnotesize $^1$Institute of Mathematics, Beijing Jiaotong
University, Beijing 100044, P. R. China
\\$^2$Department of Mathematics, Ningbo University, Ningbo 315211, P. R. China
\\$^3$Department of Computer Science, Lakehead University, Thunder Bay, Ontario, P7B 5E1, Canada}\\
  {\footnotesize
tfeng@bjtu.edu.cn, wangxiaomiao@nbu.edu.cn, rwei@lakeheadu.ca}
\vskip12pt

\end{center}

\vskip12pt

\noindent {\bf Abstract:} We consider the existence problem for a semi-cyclic holey group divisible design of type $(n,m^t)$ with block size $3$, which is denoted by a $3$-SCHGDD of type $(n,m^t)$. When $t$ is odd and $n\neq 8$ or $t$ is doubly even and $t\neq 8$, the existence problem is completely solved; when $t$ is singly even, many infinite families are obtained. Applications of our results to two-dimensional balanced sampling plans and optimal two-dimensional optical orthogonal codes are also discussed.

\noindent {\bf Keywords}: holey group divisible design; semi-cyclic; two-dimensional balanced sampling plan; two-dimensional optical orthogonal code

%%%%%%%%%%%%%%%%%%%%%%%%%%%%%%%%%%%%%%%%%%%%%%%%%%%%%%%%%%%%%%%%%%%%%%%%

\section{Introduction}

Group divisible designs are one of the most important combinatorial structures, which were widely used in constructing other combinatorial configurations \cite{bjl}. Let $K$ be a set of positive integers. A {\em group divisible
design} (GDD), denoted $K$-GDD, is a triple ($X, {\cal G},{\cal B}$)
satisfying the following properties:
\begin{enumerate}
\item[(1)] $X$ is a finite set of {\em points};

\item[(2)] $\cal G$ is a partition of $X$ into subsets (called {\em groups});

\item[(3)] $\cal B$ is a set of subsets (called {\em blocks}) of $X$, each of
cardinality from $K$, such that every $2$-subset of $X$ is either
contained in exactly one block or in exactly one group, but not in
both.
\end{enumerate}
If $\cal G$ contains $u_i$ groups of size $g_i$ for $1\leq
i\leq r$, then we call $g_1^{u_1}g_2^{u_2}\cdots g_r^{u_r}$ the {\em
group type} (or {\em type}) of the GDD. If $K=\{k\}$, we write a
$\{k\}$-GDD as a $k$-GDD. A $K$-$GDD$ of type $1^v$ is commonly called a
{\em pairwise balanced design}, denoted by a $(v,K,1)$-$PBD$. When $K=\{k\}$, a pairwise balanced design is called a {\em balanced incomplete block design}, denoted by a $(v,
k,1)$-BIBD.

A sub-GDD $(Y,{\cal H},{\cal A})$ of a GDD $(X,{\cal G},{\cal B})$ is a GDD satisfying that $Y\subseteq X$, ${\cal A}\subseteq{\cal B}$, and every group of $\cal H$ is contained in some group of $\cal G$. If a GDD has a missing sub-GDD, then we say that the GDD has a {\em hole}. In fact, the missing sub-GDD need not exist. If a GDD has several equal-sized holes which partition the point set of the GDD, then we call it a holey GDD, or HGDD. We give a formal definition of an HGDD as follows.

Let $n,m$ and $t$ be positive integers. Let $K$ be a set of positive integers. A {\em holey group divisible design} (HGDD) $K$-HGDD of type $(n,m^t)$ is a quadruple $(X,{\cal G},{\cal H},{\cal B})$ which satisfies the following properties:
\begin{enumerate}
    \item[(1)] $X$ is a finite set of $nmt$ {\em points};
    \item[(2)] ${\cal G}$ is a partition of $X$ into $n$ subsets, (called {\em groups}), each of size $mt$;
    \item[(3)] ${\cal H}$ is another partition of $X$ into $t$ subsets, (called {\em holes}), each of size $nm$ such that $|H\cap G|=m$ for each $H\in {\cal H}$ and $G\in {\cal G}$;
    \item[(4)] $\cal B$ is a set of subsets (called {\em blocks}) of $X$, each of
    cardinality from $K$, such that no block contains two distinct points of any group or any hole, but any other pair of distinct points of $X$ occurs in exactly one block of $\cal B$.
\end{enumerate}
When $m=1$, a $K$-HGDD of type $(n,1^t)$ is often said to be a {\em modified group divisible design}, denoted by a $K$-MGDD of type $t^n$. If $K=\{k\}$, we write a $\{k\}$-HGDD as a $k$-HGDD, and a $\{k\}$-MGDD as a $k$-MGDD.

Assaf \cite{a} first introduced the notion of MGDDs and settled the existence of $3$-MGDDs. The existence of $4$-MGDDs was investigated in \cite{aw,gww,lc}. There are also some results on $5$-MGDDs in \cite{aa}. Wei \cite{wei} first introduced the concept of HGDDs and gave a complete existence theorem for $3$-HGDDs. The existence of $4$-HGDDs has been completely settled in \cite{cww,gw}. We only quote the following result for the later use.

\begin{Theorem}\label{3-HGDD}{\rm \cite{wei}}
There exists a $3$-HGDD of type $(n,m^t)$ if and only if $n,t\geq 3$, $(t-1)(n-1)m\equiv 0\ ({\rm mod}\ 2)$ and $t(t-1)n(n-1)m^2\equiv 0\ ({\rm mod}\ 3)$.
\end{Theorem}

HGDDs can be seen as a special case of double group divisible designs (DGDDs), which are introduced in \cite{zhu} to simplify Stinson's proof \cite{s} on a recursive construction for group divisible designs. HGDDs can also be considered as a generalization of holey mutually orthogonal Latin squares (HMOLS), since a $k$-HGDD of type $(k,m^t)$ is equivalent to $k-2$ HMOLSs of type $m^t$ (see \cite{sz,wei} for details). HGDDs have been used for various types of combinatorial objects, such as covering designs and packing designs \cite{a2}, balanced sampling designs \cite{bcrw,gww}, etc.

In what follows we always assume that $I_u=\{0,1,\ldots,u-1\}$ and denote by $Z_v$ the additive group of integers modulo $v$. A way to construct $k$-HGDDs of type $(n,m^t)$ is the pure and mixed difference method. Let $S=\{0,t,\ldots,(m-1)t\}$ be a subgroup of order $m$ in $Z_{mt}$, and $S_l=S+l$ be a coset of $S$ in $Z_{mt}$, $0\leq l\leq t-1$. Let $X=I_n\times Z_{mt}$, ${\cal G}=\{\{i\}\times Z_{mt}:i\in I_n\}$, and ${\cal H}=\{I_n\times S_l:0\leq l\leq t-1\}$. Take a family ${\cal B}^*$ of some $k$-subsets (called {\em base blocks}) of $X$. For $i,j\in I_n$ and $B\in{\cal B}^*$, define a multi-set $\Delta_{ij}(B)=\{x-y\ ({\rm mod}\ mt): (i,x),(j,y)\in B, (i,x)\neq(j,y)\}$, and a multi-set $\Delta_{ij}({\cal B}^*)=\bigcup_{B\in{\cal B}^*}\Delta_{ij}(B)$. If for $(i,j)\in I_n\times I_n$,
$$\Delta_{ij}({\cal B}^*)=\left\{
\begin{array}{lll}
Z_{mt}\setminus S, & \ \ i\neq j,\\
\emptyset, & \ \ i=j,\\
\end{array}
\right.
$$
then a $k$-HGDD of type $(n,m^t)$ with the point set $X$, the group set $\cal G$ and the hole set $\cal H$ can be generated from ${\cal B}^*$. The required blocks are
obtained by developing all base blocks of ${\cal B}^*$ by successively
adding $1$ to the second component of each point of these base
blocks modulo $mt$. Usually a $k$-HGDD obtained by this
manner is said to be a {\em semi-cyclic $k$-HGDD} and denoted by a $k$-SCHGDD.

In this paper we shall focus on the existence of $3$-SCHGDDs of type $(n,m^t)$ with $t$ even. When $t$ is odd, the existence problem has been discussed in \cite{fwc} recently.

\begin{Lemma}\label{odd main results}{\rm \cite{fwc}}
\begin{enumerate}
\item[$(1)$] There exists a $3$-SCHGDD of type $(3,m^t)$ if and only if $(t-1)m\equiv 0\ ({\rm mod}\ 2)$ and $t\geq 3$ with the exception of $m\equiv 0\ ({\rm mod}\ 2)$ and $t=3$.
\item[$(2)$] Assume that $t\equiv 1\ ({\rm mod}\ 2)$ and $n\geq 4$. There exists a $3$-SCHGDD of type $(n,m^t)$ if and only if $t\geq 3$ and $(t-1)n(n-1)m\equiv 0\ ({\rm mod}\ 6)$ except when $(n,m,t)=(6,1,3)$, and possibly when $(1)$ $n=6$, $m\equiv 1,5\ ({\rm mod}\ 6)$ and $t\equiv 3,15\ ({\rm mod}\ 18)$, $(2)$ $n=8$, $m\equiv 2,10\ ({\rm mod}\ 12)$ and $t\equiv 7\ ({\rm mod}\ 12)$.
\end{enumerate}
\end{Lemma}

Simple counting shows that the number of base blocks in a $3$-SCHGDD of type $(n,m^t)$ is $(t-1)n(n-1)m/6$. Hence combining the result of Theorem \ref{3-HGDD}, we have the following necessary condition for the existence of $3$-SCHGDDs.

\begin{Lemma}\label{nece}
If there exists a $3$-SCHGDD of type $(n,m^t)$, then $n,t\geq 3$, $(t-1)(n-1)m\equiv 0\ ({\rm mod}\ 2)$ and $(t-1)n(n-1)m\equiv 0\ ({\rm mod}\ 6)$.
\end{Lemma}

We are to prove in Sections $3$ and $4$ that when $t$ is odd and $n\neq 8$ or $t$ is doubly even and $t\neq 8$, the existence problem for $3$-SCHGDDs is completely solved, and when $t$ is singly even, many infinite families are obtained.

Although the existence problem for SCHGDDs is interesting in its own right, there are actually some nice applications such as in the construction of two-dimensional balanced sampling plans and optimal two-dimensional optical orthogonal codes. We shall consider these applications in Section $5$. Actually, the  main result of this paper has been used for constructing optimal three dimensional optical orthogonal codes by Wang and Chang \cite{wlc}.

%To conclude an application of SCHGDDs to some interesting sampling designs are mentioned in Section $5$.

\section{Preliminaries}

Now we introduce some basic concepts and terminologies in
combinatorial design theory adopted in this paper.

An {\em automorphism} $\pi$ of a GDD $(X,{\cal G},{\cal B})$ is a permutation on $X$ leaving ${\cal G}$, ${\cal B}$ invariant, respectively. Let $H$ be the cyclic group
generated by $\pi$ under the compositions of permutations. Then all blocks of the GDD can be partitioned into some block orbits under $H$. Choose any fixed block from each block orbit and then call it a {\em base block} of this GDD under $H$. The number of the blocks contained in a block orbit is called the {\em length} of the block orbit.

A $K$-GDD of type $m^n$ is said to be {\em cyclic}, if it admits an automorphism consisting of a cycle of length $mn$. A cyclic $K$-GDD is denoted by a $K$-CGDD. For a $K$-CGDD $(X,{\cal G},{\cal B})$, we can always identify $X$ with $Z_{mn}$ and ${\cal G}$ with $\{\{in+j : 0\leq i\leq m-1\} : 0\leq j\leq n-1\}$. If the length of each block orbit in a $K$-CGDD of type $m^n$ is $mn$, then the $K$-GDD is called {\em strictly cyclic}.

\begin{Lemma}\label{3-sCGDD}{\rm \cite{wc}}
There exists a strictly cyclic $3$-GDD of type $m^n$ if and only if
\begin{enumerate}
\item[$(1)$] $m(n-1)\equiv 0\ ({\rm mod}\ 6)$ and $n\geq 4;$
\item[$(2)$] $n\not\equiv 2,3\ ({\rm mod}\ 4)$ when $m\equiv 2\ ({\rm mod}\ 4).$
\end{enumerate}
\end{Lemma}

\begin{Construction} \label{SCHGDD-from strictly CGDD}{\rm \cite{fwc}}
Suppose that there exist a strictly cyclic $K$-GDD of type $w^t$ and an $l$-MGDD of type $k^n$ for each $k\in K$. Then there exists an $l$-SCHGDD of type $(n,w^t)$.
\end{Construction}

Construction \ref{SCHGDD-from strictly CGDD} shows that strictly cyclic $K$-GDDs are helpful to yield SCHGDDs. Combining the results of Corollary $4.10$ and Lemma $4.15$ in \cite{fwc}, we have the following strictly cyclic GDDs.

\begin{Lemma}\label{4^t-0,1,2mod3} {\rm \cite{fwc}}
There exists a strictly cyclic $\{3,5\}$-GDD of type $4^t$ for any $t\geq 4$ and $t\not\in\{5,8,11\}$.
\end{Lemma}

Let $K$ be a set of positive integers and $v$ be an odd positive integer. For each $1\leq i\leq b$, let $B_i\subseteq Z_v$ and $|B_i|=k_i\in K$. The $b$ sets $B_i=\{x_{i,1}, x_{i,2},\ldots,x_{i,k}\}$, $1\leq i\leq b$, form a {\em perfect $(v,K, 1)$ difference family}, written as a $(v,K,1)$-PDF, if all the differences $x_{i,s}-x_{i,r}$, $1\leq i\leq b$, $1\leq r < s \leq k_i$ cover exactly the set $\{1,2,\ldots,(v-1)/2\}$.

\begin{Lemma}\label{relation PDF and GDD} {\rm \cite{fwc}} If there exists a $(2t-1,K,1)$-PDF, then there exists a strictly cyclic $K$-GDD of type $2^t$.
\end{Lemma}

\begin{Lemma}
{\rm \cite{czl} }\label{PDF-24-19}
There exists a $(v,\{3,4\},1)$-PDF for $v\equiv 1\ ({\rm mod}\ 6)$ and $v\geq 19$.
\end{Lemma}

A $K$-GDD of type $m^n$ is said to be {\em semi-cyclic}, if it admits an automorphism which permutes the elements of each group $G\in{\cal G}$ in an $m$ cycle. Such a GDD is denoted by a $K$-SCGDD of type $m^n$. For a $K$-SCGDD $(X,{\cal G},{\cal B})$, we can always identify $X$ with $I_n\times Z_m$ and ${\cal G}$ with $\{\{i\}\times Z_m:i\in I_n\}$. In this case the automorphism can be taken as $(i,x)\longmapsto(i,x+1)$ (mod $(-,m)$), $i\in I_n$ and $x\in Z_m$. Assume that ${\cal B}^*$ is the set of base blocks of a $K$-SCGDD of type $m^n$. It is easy to verify that for $(i,j)\in I_n\times I_n$,
$$\Delta_{ij}({\cal B}^*)=\left\{
\begin{array}{lll}
Z_m, & \ \ \ i\neq j,\\
\emptyset, & \ \ \ i=j.\\
\end{array}
\right.
$$
Note that no element of $Z_m$ occurs more than once in $\Delta_{ij}({\cal B}^*)$, and the length of each block orbit is $m$ in a $K$-SCGDD.

\begin{Lemma}{\rm\cite{jiang}}\label{scgdd}
There is a $3$-SCGDD of type $m^n$ if and only if $n\geq 3$ and
\begin{enumerate}
\item[$(1)$] $(n-1)m\equiv 0\ ({\rm mod}\ 2)$;
\item[$(2)$] $n(n-1)m\equiv 0\ ({\rm mod}\ 3)$;
\item[$(3)$] $n\not\equiv 2,3\ ({\rm mod}\ 4)$ when $m\equiv 2\ ({\rm mod}\ 4)$, and $n\neq 3$ when $m\equiv 0\ ({\rm mod}\ 4)$.
\end{enumerate}
\end{Lemma}

\begin{Construction} \label{SCHGDD-recur}{\rm \cite{fwc}}
Suppose that there exist a $K$-SCGDD of type $g^n$ and an $l$-SCHGDD of type $(k,w^t)$ for each $k\in K$. Then there exists an $l$-SCHGDD of type $(n,(gw)^t)$.
\end{Construction}

A $(k,m)$-CDM is a $k\times m$ matrix $D=(d_{ij})$ with entries from $Z_m$ such that for any two distinct rows $x$ and $y$, the difference list $\{d_{xj}-d_{yj}:j\in I_m\}$ contains each integer of $Z_m$ exactly once. It is  known that a $(3,m)$-CDM exists for any positive odd integer $m$ \cite{bjl}.

\begin{Construction}\label{SCHGDD from CDM}{\rm \cite{fwc}}
If there exist a $k$-SCHGDD of type
$(n,w^t)$ and a $(k,v)$-CDM, then there exists a $k$-SCHGDD of type $(n,(wv)^t)$.
\end{Construction}

The following construction is simple but useful.

\begin{Construction}\label{SCHGDD from SCHGDD}{\rm \cite{fwc}}
If there exist a $k$-SCHGDD of type $(n,(gw)^t)$ and a $k$-SCHGDD of type
$(n,g^w)$, then there exists a $k$-SCHGDD of type $(n,g^{wt})$.
\end{Construction}

As the end of this section, we state a shorter definition of SCHGDDs. An {\em automorphism} of an HGDD $(X,{\cal G},{\cal H},{\cal B})$ is a permutation on $X$ leaving ${\cal G}$, ${\cal H}$, $\cal B$ invariant, respectively. A $K$-SCHGDD of type
$(n,m^t)$ is a $K$-HGDD of type $(n,m^t)$
with an automorphism of order $mt$ that permutes the points within each group.

\section{$3$-SCHGDDs of type $(n,m^t)$ for $n=4,5,6,8$}

In this section, we always assume that $[a,b]$ denotes the set of integers $n$ such that $a\leq n\leq b$, and $[a,b]_{\not\equiv r(v)}$, $0\leq r\leq v-1$, denotes the set of integers $n$ such that $a\leq n\leq b$ and $n\not\equiv r\ ({\rm mod}\ v)$.

\begin{Lemma}\label{4-SCHGDD}
If $(t-1)m\equiv 0\ ({\rm mod}\ 2)$ and $t\geq 3$, then there exists a $3$-SCHGDD of type $(4,m^t)$.
\end{Lemma}

\proof When $t$ is odd, the conclusion follows from Lemma \ref{odd main results}. When $m\equiv 0\ ({\rm mod}\ 4)$, start from a $3$-SCGDD of type $(m/2)^4$, which exists by Lemma \ref{scgdd}. By Lemma \ref{odd main results}, there exists a $3$-SCHGDD of type $(3,2^t)$ for any even integer $t\geq 4$. Then apply Construction \ref{SCHGDD-recur} to obtain a $3$-SCHGDD of type $(4,m^t)$.

When $m=2$ and $t\equiv 0\ ({\rm mod}\ 2)$, the $4(t-1)$ base blocks of a $3$-SCHGDD of type $(4,2^t)$ on $I_4\times Z_{2t}$ with the group set $\{\{i\}\times Z_{2t}:i\in I_4\}$ and the hole set $\{I_4\times \{j,t+j\}:0\leq j\leq t-1\}$ are listed as follows:

\begin{center}
\begin{tabular}{l}
$\{(0,0),(1,i),(2,2i)\}$, $i\in[1,t-1]\setminus\{t/2-1,t/2\}$;\\
$\{(0,0),(1,t+i),(3,t-i)\}$, $i\in[1,t-1]\setminus\{t/2,t/2+1\}$;\\
$\{(0,0),(2,2i+1),(3,t+i)\}$, $i\in[1,t-2]$;\\
$\{(1,0),(2,t+i),(3,2i+1)\}$, $i\in[1,t-2]$;\\
\end{tabular}

\begin{tabular}{ll}
$\{(0,0),(1,3t/2+1),(2,1)\},$ & $\{(0,0),(1,t/2-1),(3,t/2)\},$ \\ $\{(0,0),(2,t-2),(3,2t-1)\},$&
$\{(0,0),(1,3t/2),(2,2t-1)\},$ \\ $\{(0,0),(1,t/2),(3,t/2-1)\},$ & $\{(1,0),(2,2t-1),(3,t-2)\}.$\\
\end{tabular}
\end{center}

\noindent When $m\equiv 2\ ({\rm mod}\ 4)$ and $t\equiv 0\ ({\rm mod}\ 2)$, start from a $3$-SCHGDD of type $(4,2^t)$, and apply Construction \ref{SCHGDD from CDM} with a $(3,m/2)$-CDM to obtain a $3$-SCHGDD of type $(4,m^t)$. \qed

\begin{Lemma}\label{5-SCHGDD-1^4}
There is no $3$-SCHGDD of type $(5,1^4)$.
\end{Lemma}

\proof Suppose that there exists a $3$-SCHGDD of type $(5,1^4)$ on $I_5\times Z_4$ with the group set $\{\{(i,0),(i,1),(i,2),(i,3)\}:0\leq i\leq 4\}$ and the hole set $\{\{(i,j):0\leq i\leq 4\}:0\leq j\leq 3\}$. It has $10$ base blocks. Denote by ${\cal B}^*$ the set of its base blocks. Let ${\cal B}^*=\{\{(a_1^{(l)},z_1^{(l)}),(a_2^{(l)},z_2^{(l)}),(a_3^{(l)},z_3^{(l)})\}:0\leq l\leq 9\}$. Given $b,c\in I_5$ and $b\neq c$, count the number of base blocks containing the pairs of the form $\{(b,*),(c,*)\}$. The number is $3$. Now write ${\cal A}=\{\{a_1^{(l)},a_2^{(l)},a_3^{(l)}\}:0\leq l\leq 9\}$. Then ${\cal A}$ forms a $(5,3,3)$-BIBD on $I_5$. Up to isomorphism, one can easily check that there is only one $(5,3,3)$-BIBD (hint: each point in a $(5,3,3)$-BIBD occurs in exactly $6$ blocks). Hence, without loss of generality the $10$ base blocks are assumed as follows:
\begin{center}\begin{tabular}{lll}
$\{(0,0),(1,x_0),(2,y_0)\}$,&
$\{(0,0),(1,x_1),(3,y_1)\}$,&
$\{(0,0),(1,x_2),(4,y_2)\}$,\\
$\{(0,0),(2,x_3),(3,y_3)\}$,&
$\{(0,0),(2,x_4),(4,y_4)\}$,&
$\{(0,0),(3,x_5),(4,y_5)\}$,\\
$\{(1,0),(2,x_6),(3,y_6)\}$,&
$\{(1,0),(2,x_7),(4,y_7)\}$,&
$\{(1,0),(3,x_8),(4,y_8)\}$,\\
$\{(2,0),(3,x_9),(4,y_9)\}$,
\end{tabular}
\end{center}
where for any $0\leq l\leq 9$, $(x_l,y_l)\in\{(1,2),(2,1),(1,3),(3,1),(2,3),(3,2)\}$. Let $i_1,i_2\in I_5$ and $i_1\neq i_2$. Write $\Delta_{i_1 i_2}=\{x-y\ ({\rm mod}\ 4):\{(i_1,x),(i_2,y)\}\subset B, B\in{\cal B}^*\}$. Then we have
\begin{center}
\begin{tabular}{lll}
$\Delta_{10}=\{x_0,x_1,x_2\}$,&
$\Delta_{21}=\{x_6,x_7,y_0-x_0\}$,&
$\Delta_{32}=\{x_9,y_3-x_3,y_6-x_6\}$,\\
$\Delta_{20}=\{x_3,x_4,y_0\}$,&
$\Delta_{31}=\{x_8,y_6,y_1-x_1\}$,&
$\Delta_{42}=\{y_9,y_4-x_4,y_7-x_7\}$,\\
$\Delta_{30}=\{x_5,y_1,y_3\}$,&
$\Delta_{41}=\{y_7,y_8,y_2-x_2\}$,&
$\Delta_{43}=\{y_5-x_5,y_8-x_8,y_9-x_9\}$,\\
$\Delta_{40}=\{y_2,y_4,y_5\}$.
\end{tabular}
\end{center}
Without loss of generality, assume that $(x_0,x_1,x_2)=(1,2,3)$. It follows that $(y_0,y_1,y_2)\in\{(2,1,1),(2,1,2),(2,3,1),(2,3,2),(3,1,1),(3,1,2),(3,3,1),(3,3,2)\}$. For convenience we make the following table,
\begin{center}
{\small
\begin{tabular}{|c|c|c|c|c|}
\hline
$(y_0,y_1,y_2)$   & $(x_3,x_4,x_5,y_3,y_4,y_5)$ & $(y_3-x_3,y_4-$ & $(x_6,x_7,x_8,y_6,y_7,y_8)$ & $(y_6-x_6,y_7-$\\
     & $$& $x_4,y_5-x_5)$ & $$& $x_7,y_8-x_8)$\\
\hline
 $(2,1,1)$  & $(1,3,2,3,2,3)$ & $(2,3,1)$ & $(2,3,2,1,1,3)$ & $(3,2,1)$\\

  & $(3,1,3,2,3,2)$ & $(3,2,3)$ & $(3,2,1,2,1,3)$ & $(3,3,2)$\\

 &
 &  & $(3,2,2,1,1,3)$ & $(2,3,1)$\\

  &  & & $(3,2,2,1,3,1)$ & $(2,1,3)$\\
\hline
 $(2,1,2)$  & $(1,3,2,3,1,3)$ & $(2,2,1)$ & $(2,3,2,1,2,1)$ & $(3,3,3)$\\

 & $(3,1,3,2,3,1)$
 & $(3,2,2)$ & $(3,2,1,2,1,2)$ & $(3,3,1)$\\
\hline
 $(2,3,1)$  & $(1,3,1,2,2,3)$ & $(1,3,2)$ & $(2,3,2,3,1,3)$&$(1,2,1)$\\

  & $(3,1,1,2,3,2)$ & $(3,2,1)$ & $(3,2,3,2,3,1)$ &$(3,1,2)$ \\

 & $(3,1,1,2,2,3)$
 & $(3,1,2)$ &  & \\

  & $(3,1,2,1,2,3)$ & $(2,1,1)$ & &\\
\hline
 $(2,3,2)$  & $(1,3,1,2,1,3)$ & $(1,2,2)$ &$(2,3,2,3,2,1)$ & $(1,3,3)$\\

 & $(3,1,2,1,3,1)$
 &$(2,2,3)$& $(3,2,3,2,1,2)$ & $(3,3,3)$\\
\hline
 $(3,1,1)$  & $(1,2,3,2,3,2)$ & $(1,1,3)$& $(1,3,1,2,1,3)$ & $(1,2,2)$\\

 & $(2,1,2,3,2,3)$
 & $(1,1,1)$ &$(3,1,2,1,3,1)$ & $(2,2,3)$\\
\hline
 $(3,1,2)$  & $(1,2,2,3,1,3)$ & $(2,3,1)$ & $(1,3,1,2,1,2)$&$(1,2,1)$\\

  & $(1,2,2,3,3,1)$ & $(2,1,3)$ & $(3,1,2,1,2,1)$ &$(2,1,3)$ \\

 & $(1,2,3,2,3,1)$
 & $(1,1,2)$&  & \\

  & $(2,1,2,3,3,1)$ & $(1,2,3)$ & &\\
\hline
 $(3,3,1)$  & $(1,2,1,2,3,2)$ & $(1,1,1)$& $(1,3,2,3,1,3)$ & $(2,2,1)$\\

 & $(2,1,2,1,2,3)$
 & $(3,1,1)$& $(3,1,3,2,3,1)$ & $(3,2,2)$\\
\hline
 $(3,3,2)$  & $(1,2,1,2,1,3)$ &$(1,3,2)$ & $(1,3,2,3,2,1)$ & $(2,3,3)$\\

  & $(2,1,2,1,3,1)$ & $(3,2,3)$ & $(1,3,3,2,1,2)$& $(1,2,3)$\\

 &
 & & $(1,3,3,2,2,1)$ & $(1,3,2)$\\

  &  && $(3,1,3,2,2,1)$ & $(3,1,2)$\\
\hline
\end{tabular}
}\end{center}
where the values of $x_3,x_4,x_5,y_3,y_4,y_5$ in the $2^{{\rm th}}$ column are determined by $\Delta_{20}=\Delta_{30}=\Delta_{40}=\{1,2,3\}$, and the values of $x_6,x_7,x_8,y_6,y_7,y_8$ in the $4^{{\rm th}}$ column are determined by $\Delta_{21}=\Delta_{31}=\Delta_{41}=\{1,2,3\}$. Then analyze the $3^{{\rm th}}$ and the $5^{{\rm th}}$ columns with the condition $\Delta_{32}=\Delta_{42}=\Delta_{43}=\{1,2,3\}$. We can see that no $x_9,y_9$ exist, a contradiction. \qed

\begin{Lemma}\label{5-SCHGDD-m^4}
There exists a $3$-SCHGDD of type $(5,m^4)$ for any integer $m\equiv 1,5\ ({\rm mod}\ 6)$ and $m\geq 5$.
\end{Lemma}

\proof We here construct a $3$-SCHGDD of type $(5,m^4)$ on $Z_5\times Z_{4m}$ with the group set $\{\{i\}\times Z_{4m}:i\in Z_5\}$ and the hole set $\{Z_5\times \{j,4+j,\ldots,4(m-1)+j\}:0\leq j\leq 3\}$. All the $10m$ base blocks can be obtained from the following $2m$ initial blocks by $(+1\ {\rm mod}\ 5,-)$.

When $m=5$, the $10$ initial blocks are:
\begin{center}
\begin{tabular}{lll}
$\{(0,0),(1,6),(2,19)\}$, & $\{(0,0),(1,1),(2,18)\}$, &
$\{(0,0),(1,3),(2,17)\}$, \\ $\{(0,0),(1,10),(2,9)\}$, &
$\{(0,0),(1,9),(2,7)\}$, & $\{(0,0),(1,2),(3,7)\}$, \\
$\{(0,0),(1,5),(3,6)\}$, & $\{(0,0),(1,7),(3,9)\}$, &
$\{(0,0),(1,11),(3,14)\}$, \\ $\{(0,0),(1,15),(3,5)\}$. &&\\
\end{tabular}
\end{center}
When $m\equiv 1,5\ ({\rm mod}\ 6)$ and $m\geq7$, the $2m$ initial blocks are divided into two parts. The first part consists of the following $m+2$ initial blocks:

\begin{center}\begin{tabular}{ll}
$\{(0,0),(1,2m-4-4i),(3,4m-3-2i)\}$, & $i\in[0,m-2]\setminus\{(m-1)/2\}$;\\
\end{tabular}

\begin{tabular}{ll}
$\{(0,0),(1,2m-1),(3,4m-1)\}$, & $\{(0,0),(1,4m-1),(3,1)\}$,\\
$\{(0,0),(1,2m),(2,3m)\}$, & $\{(0,0),(1,4m-2),(2,m+2)\}$.
\end{tabular}
\end{center}
The second part consists of the following $m-2$ initial blocks:

\begin{itemize}
%$\cdot$
\item if $m\equiv 7,11\ ({\rm mod}\ 12)$, then
\begin{center}\begin{tabular}{ll}
$\{(0,0),(1,2m+1-4i),(2,2m-4-8i)\}$,& $i\in[0,(m-7)/4]$;\\
$\{(0,0),(1,2m-5-4i),(2,2m-8-8i)\}$,& $i\in[0,(m-7)/4]$;\\
$\{(0,0),(1,1+2i),(2,2m+4+4i)\}$,& $i\in[0,(m-3)/2]$;
\end{tabular}
\end{center}

%$\cdot$

\item
if $m\equiv 1,5\ ({\rm mod}\ 12)$ and $m\geq 13$, firstly we take $(m+3)/2$ initial blocks:
\begin{center}\begin{tabular}{lll}
$\{(0,0),(1,2m+1-4i),(2,2m-4-8i)\}$, &$i\in[0,(m-9)/4]$;\\
$\{(0,0),(1,2m-5-4i),(2,2m-8-8i)\}$,& $i\in[0,(m-13)/4]$;\\
\end{tabular}

\begin{tabular}{lll}
$\{(0,0),(1,m-2),(2,2m+4)\}$, & $\{(0,0),(1,3m+2),(2,2m+8)\}$,\\
$\{(0,0),(1,1),(2,6)\}$, & $\{(0,0),(1,2m+3),(2,10)\}$,\\
$\{(0,0),(1,9),(2,2m+20)\}$;
\end{tabular}
\end{center}
\noindent and then when $m\equiv 1\ ({\rm mod}\ 12)$ and $m\geq 13$, we take the $(m-7)/2$ initial blocks:

\begin{center}\begin{tabular}{lll}
$\{(0,0),(1,15+6i),(2,2m+24+12i)\}$,& $i\in[0,(m-13)/6]$;\\
$\{(0,0),(1,13+6i),(2,2m+32+12i)\}$, &$i\in[0,(m-19)/6]$, & (null if $m=13$);\\
$\{(0,0),(1,11+6i),(2,2m+28+12i)\}$,& $i\in[0,(m-19)/6]$, & (null if $m=13$);
\end{tabular}

\begin{tabular}{lll}
$\{(0,0),(1,3),(2,2m+16)\}$, & $\{(0,0),(1,7),(2,2m+12)\}$;
\end{tabular}
\end{center}
\noindent when $m\equiv 5\ ({\rm mod}\ 12)$ and $m\geq 17$, we take the following $(m-7)/2$ initial blocks:

\begin{center}\begin{tabular}{lll}
$\{(0,0),(1,19+6i),(2,2m+32+12i)\}$,& $i\in[0,(m-17)/6]$;\\
$\{(0,0),(1,17+6i),(2,2m+40+12i)\}$, &$i\in[0,(m-23)/6]$, & (null if $m=17$);\\
$\{(0,0),(1,15+6i),(2,2m+36+12i)\}$,& $i\in[0,(m-23)/6]$, & (null if $m=17$);
\end{tabular}

\begin{tabular}{lll}
$\{(0,0),(1,3),(2,2m+12)\}$, & $\{(0,0),(1,7),(2,2m+24)\}$,\\
$\{(0,0),(1,11),(2,2m+16)\}$, & $\{(0,0),(1,13),(2,2m+28)\}$.
\end{tabular}
\end{center}
\end{itemize}

\begin{Lemma}\label{5-SCHGDD}
If $(t-1)m\equiv 0\ ({\rm mod}\ 3)$ and $t\geq 3$, then there exists a $3$-SCHGDD of type $(5,m^t)$ with the exception of $(m,t)=(1,4)$.
\end{Lemma}

\proof When $t$ is odd, the conclusion follows from Lemma \ref{odd main results}. When $m\equiv 0\ ({\rm mod}\ 6)$, start from a $3$-SCGDD of type $(m/2)^5$, which exists by Lemma \ref{scgdd}. By Lemma \ref{odd main results}, there exists a $3$-SCHGDD of type $(3,2^t)$ for any even integer $t\geq 4$. Then apply Construction \ref{SCHGDD-recur} to obtain a $3$-SCHGDD of type $(5,m^t)$.

When $m=3$ and $t\equiv 0\ ({\rm mod}\ 2)$, the $10(t-1)$ base blocks of a $3$-SCHGDD of type $(5,3^t)$ on $I_5\times Z_{3t}$ with the group set $\{\{i\}\times Z_{3t}:i\in I_5\}$ and the hole set $\{I_5\times \{j,t+j,2t+j\}:0\leq j\leq t-1\}$ are listed as follows:

\begin{center}
\begin{tabular}{l}
$\{(l,0),(l+1,i),(l+2,2i+2t)\}$, $i\in[1,t-1]\setminus\{t/2-1,t/2\}$;\\
$\{(l,0),(l+1,t+i),(l+3,2t-i+1)\}$, $i\in[2,t-1]$;\\
\end{tabular}

\begin{tabular}{ll}
$\{(l,0),(l+1,5t/2),(l+2,2t-1)\},$ & $\{(l,0),(l+1,t/2),(l+2,t-1)\},$ \\ $\{(l,0),(l+1,t+1),(l+3,2)\}$,\\
\end{tabular}
\end{center}

\noindent where $l$ runs over $I_5$, and the first components of elements are reduced modulo $5$. When $m\equiv 3\ ({\rm mod}\ 6)$ and $t\equiv 0\ ({\rm mod}\ 2)$, start from a $3$-SCHGDD of type $(5,3^t)$, and apply Construction \ref{SCHGDD from CDM} with a $(3,m/3)$-CDM to obtain a $3$-SCHGDD of type $(5,m^t)$.

When $m=2$ and $t\equiv 10\ ({\rm mod}\ 12)$, take a $(2t-1,\{3,4\},1)$-PDF from Lemma \ref{PDF-24-19}, which yields a strictly cyclic $\{3,4\}$-GDD of type $2^t$ by Lemma \ref{relation PDF and GDD}. Start from this GDD, and apply Construction \ref{SCHGDD-from strictly CGDD} with a $3$-MGDD of type $k^5$, $k\in\{3,4\}$, which exists by Theorem \ref{3-HGDD}, to obtain a $3$-SCHGDD of type $(5,2^t)$. When $m\equiv 2,10\ ({\rm mod}\ 12)$ and $t\equiv 10\ ({\rm mod}\ 12)$, start from a $3$-SCHGDD of type $(5,2^t)$, and apply Construction \ref{SCHGDD from CDM} with a $(3,m/2)$-CDM to obtain a $3$-SCHGDD of type $(5,m^t)$.

When $m\equiv 2,10\ ({\rm mod}\ 12)$ and $t\equiv 4\ ({\rm mod}\ 12)$, or $m\equiv 4,8\ ({\rm mod}\ 12)$ and $t\equiv 4\ ({\rm mod}\ 6)$, by Lemma \ref{3-sCGDD}, we have a strictly cyclic $3$-GDD of type $m^t$. Then apply Construction \ref{SCHGDD-from strictly CGDD} with a $3$-MGDD of type $3^5$ to obtain a $3$-SCHGDD of type $(5,m^t)$.

When $m=1$, $t\equiv 4\ ({\rm mod}\ 6)$ and $t\geq 10$, the $10(t-1)/3$ base blocks of a $3$-SCHGDD of type $(5,1^t)$ on $I_5\times Z_t$ with the group set $\{\{i\}\times Z_t:i\in I_5\}$ and the hole set $\{I_5\times \{j\}:0\leq j\leq t-1\}$ are listed as follows:

\begin{center}
\begin{tabular}{l}
$\{(l,0),(l+1,i),(l+2,2(t-1)/3+2i)\}$, $i\in[1,(t-4)/3]\setminus\{(t+2)/6\}$;\\
$\{(l,0),(l+1,(t-1)/3+i),(l+3,2(t-1)/3-i)\}$, $i\in[1,(t-4)/3]\setminus\{(t-4)/6\}$;\\
\end{tabular}

\begin{tabular}{ll}
$\{(l,0),(l+1,(t+2)/6),(l+2,t/2)\},$ & $\{(l,0),(l+1,t/2-1),(l+2,(t-4)/3)\},$ \\
$\{(l,0),(l+1,2(t-1)/3),(l+3,t-1)\},$ & $\{(l,0),(l+1,t-1),(l+3,(t-1)/3)\},$ \\
\end{tabular}
\end{center}

\noindent where $l$ runs over $I_5$, and the first components of elements are reduced modulo $5$. When $m\equiv 1,5\ ({\rm mod}\ 6)$, $t\equiv 4\ ({\rm mod}\ 6)$ and $t\geq 10$, start from a $3$-SCHGDD of type $(5,1^t)$, and apply Construction \ref{SCHGDD from CDM} with a $(3,m)$-CDM to obtain a $3$-SCHGDD of type $(5,m^t)$.

When $m=1$ and $t=4$, by Lemma \ref{5-SCHGDD-1^4}, there is no $3$-SCHGDD of type $(5,1^4)$. When $m\equiv 1,5\ ({\rm mod}\ 6)$, $m\geq 5$ and $t=4$, the conclusion follows from Lemma \ref{5-SCHGDD-m^4}. \qed

\begin{Lemma}
\label{6,(2,4)^8} There exists a $3$-SCHGDD of type $(6,m^8)$ for $m=2,4$.
\end{Lemma}

\proof We here construct a $3$-SCHGDD of type $(6,m^8)$ on $(Z_5\cup\{\infty\})\times Z_{8m}$ with the group set $\{\{i\}\times Z_{8m}:i\in Z_5\cup\{\infty\}\}$ and the hole set $\{(Z_5\cup\{\infty\})\times \{j,8+j,\ldots,8(m-1)+j\}:0\leq j\leq 7\}$. All the $35m$ base blocks can be obtained from the following $7m$ initial blocks by $(+1\ {\rm mod}\ 5,-)$, where $\infty+1=\infty$.

\begin{center}
\begin{tabular}{lll}
$m=2:$\\
$\{(0,0), (1,1), (3,2) \}$,& $\{(0,0),(3,1),(2,2)\}$,&  $\{(0,0), (\infty,1), (1,2)\}$,\\
$\{(0,0), (4,2), (\infty,4)\}$,&$\{(0,0),(1,3), (3,6)\}$,&$\{(0,0), (3,3),(2,6)\}$,\\
$\{(0,0, (\infty,3), (1,7) \}$,&$\{(0,0,(1,4), (2,9) \}$,&$\{(0,0), (2,4),(\infty,10)\}$,\\
$\{(0,0), (3,4),(1,9) \}$,& $\{(0,0),(4,4), (3,9)\}$,&  $\{(0,0), (2,5),(\infty, 14)\}$,\\
$\{(0,0), (\infty,5), (1,10)\}$,&     $\{(0,0),(1,6), (\infty,13)\}$.\\
\end{tabular}
\end{center}

\begin{center}
\begin{tabular}{lll}
$m=4:$\\
$\{(0,0),(3,4),(\infty,9)\}$,&
$\{(0,0),(4,4),(2,9)\}$,&
$\{(0,0),(1,6),(\infty,12)\}$,\\
$\{(0,0),(4,6),(1,13)\}$,&
$\{(0,0),(1,7),(4,14)\}$,&
$\{(0,0),(4,7),(1,17)\}$,\\
$\{(0,0),(\infty,7),(1,14)\}$,&
$\{(0,0),(4,9),(1,20)\}$,&
$\{(0,0),(1,10),(2,21)\}$,\\
$\{(0,0),(3,10),(2,20)\}$,&
$\{(0,0),(\infty,10),(2,6)\}$,&
$\{(0,0),(1,1),(3,2)\}$,\\
$\{(0,0),(3,1),(2,2)\}$,&
$\{(0,0),(\infty,1),(1,2)\}$,&
$\{(0,0),(4,2),(1,5)\}$,\\
$\{(0,0),(\infty,2),(1,4)\}$,&
$\{(0,0),(1,3),(2,12)\}$,&
$\{(0,0),(3,3),(1,12)\}$,\\
$\{(0,0),(4,3),(2,17)\}$,&
$\{(0,0),(\infty,3),(2,14)\}$,&
$\{(0,0),(2,4),(\infty,23)\}$,\\
$\{(0,0),(2,5),(\infty,20)\}$,&
$\{(0,0),(4,5),(\infty,22)\}$,&
$\{(0,0),(3,6),(1,19)\}$,\\
$\{(0,0),(4,11),(\infty,29)\}$,&
$\{(0,0),(\infty,11),(3,17)\}$,&
$\{(0,0),(2,13),(\infty,27)\}$,\\
$\{(0,0),(\infty,13),(4,17))\}$.
\end{tabular}
\end{center}

\begin{Lemma}\label{6,2-4^power} There exists a $3$-SCHGDD of type $(6,2^{2^r})$ for $r\geq 2$.
\end{Lemma}

\proof We use induction on $r$. When $r=2$, take a strictly cyclic $3$-GDD of type $2^4$ from Lemma \ref{3-sCGDD}. Then apply Construction \ref{SCHGDD-from strictly CGDD} with a $3$-MGDD of type $3^6$, which exists by Theorem \ref{3-HGDD}, to obtain a $3$-SCHGDD of type $(6,2^4)$. When $r=3$, the conclusion follows from Lemma \ref{6,(2,4)^8}.
When $r\geq 4$, assume that there exists a $3$-SCHGDD of type $(6,2^{2^{r-2}})$. Due to $2^x\equiv 4,8\ ({\rm mod}\ 12)$ for any $x\geq 2$, by Lemma \ref{3-sCGDD} we have a strictly cyclic $3$-GDD of type $(2^x)^4$. Then apply Construction \ref{SCHGDD-from strictly CGDD} with a $3$-MGDD of type $3^6$ to obtain a $3$-SCHGDD of type $(6,(2^x)^4)$. Start from a $3$-SCHGDD of type $(6,(2^{r-1})^4)$, and apply Construction \ref{SCHGDD from SCHGDD} with the given $3$-SCHGDD of type $(6,2^{2^{r-2}})$ to obtain a $3$-SCHGDD of type $(6,2^{2^r})$. \qed

\begin{Lemma}\label{6-SCHGDD-m^3}
There exists a $3$-SCHGDD of type $(6,m^3)$ for any integer $m\equiv 1\ ({\rm mod}\ 2)$ and $m\geq 3$.
\end{Lemma}

\proof Let $I=Z_5\cup\{\infty\}$. We here construct a $3$-SCHGDD of type $(6,m^3)$ on $I\times Z_{3m}$ with the group set $\{\{i\}\times Z_{3m}:i\in I\}$ and the hole set $\{I\times \{j,3+j,\ldots,3(m-1)+j\}:0\leq j\leq 2\}$. Only initial base blocks are listed below, and all other base blocks are obtained by developing these base blocks by $(+1,-)$ modulo $(5,-)$, where $\infty+1=\infty$. When $m\equiv 1\ ({\rm mod}\ 4)$ and $m\geq5$:
\begin{center}\begin{tabular}{lll}
$\{(0,0),(1,1+i),(2,(3m+7)/2+2i)\}$, & $i\in[0,(3m-7)/2]_{\not\equiv 2(3)}$;\\
$\{(0,0),(2,5+2i),(\infty,(3m+5)/2+i)\}$, & $i\in[0,(3m-23)/4]_{\not\equiv 2(3)}$,\\
& (null if $m=5$);\\
$\{(0,0),(2,(3m+1)/2+2i),(\infty,(9m-5)/4+i)\}$, &  $i\in[0,(3m-11)/4]_{\not\equiv 2(3)}$;
\end{tabular}
\begin{tabular}{lll}
$\{(0,0),(2,(3m-5)/2),(\infty,3m-2)\}$, & $\{(0,0),(2,3m-1),(\infty,(3m+1)/4)\}$,\\
$\{(0,0),(1,(3m+1)/2),(\infty,(3m-1)/2)\}$, & $\{(0,0),(1,(3m-1)/2),(3,(3m+1)/2)\}$.
\end{tabular}
\end{center}
When $m\equiv 3\ ({\rm mod}\ 4)$ and $m\geq7$:
\begin{center}\begin{tabular}{lll}
$\{(0,0),(1,1+i),(2,(3m+7)/2+2i)\}$, & $i\in[0,(3m-9)/2]_{\not\equiv 2(3)}$;\\
$\{(0,0),(2,4+2i),(\infty,(3m+7)/2+i)\}$, & $i\in[0,(3m-9)/4]_{\not\equiv 1(3)}$;\\
$\{(0,0),(2,(3m+13)/2+2i),(\infty,(9m+13)/4+i)\}$, & $i\in[0,(3m-21)/4]_{\not\equiv 2(3)}$;
\end{tabular}
\begin{tabular}{lll}
$\{(0,0),(2,(3m-7)/2),(\infty,3m-1)\}$, & $\{(0,0),(1,3m-1),(\infty,(3m-5)/4)\}$,\\
$\{(0,0),(1,(3m+1)/2),(\infty,1)\}$, & $\{(0,0),(1,(3m-1)/2),(3,2)\}$,\\
$\{(0,0),(1,(3m-5)/2),(3,3m-2)\}$. &
\end{tabular}
\end{center}
When $m=3$, the conclusion follows from Lemma \ref{odd main results}. \qed

\begin{Lemma}\label{6-SCHGDD-1-3^r} There exists a $3$-SCHGDD of type $(6,1^{3^r})$ for any integer $r\geq 2$.
\end{Lemma}

\proof We use induction on $r$. When $r=2$, the conclusion follows from Lemma \ref{odd main results}. When $r\geq 3$, assume that there exists a $3$-SCHGDD of type $(6,1^{3^{r-1}})$. By Lemma \ref{6-SCHGDD-m^3} we have a $3$-SCHGDD of type $(6,(3^{r-1})^3)$. Then apply Construction \ref{SCHGDD from SCHGDD} with the given $3$-SCHGDD of type $(6,1^{3^{r-1}})$, we have the required $3$-SCHGDD of type $(6,1^{3^r})$. \qed

\begin{Lemma}\label{6-SCHGDD}
If $(t-1)m\equiv 0\ ({\rm mod}\ 2)$ and $t\geq 3$, then there exists a $3$-SCHGDD of type $(6,m^t)$ with the exception of $(m,t)=(1,3)$.
\end{Lemma}

\proof When $t$ is odd with the exception of $m\equiv 1,5\ ({\rm mod}\ 6)$ and $t\equiv 3,15\ ({\rm mod}\ 18)$, the conclusion follows from Lemma \ref{odd main results}. When $m=1$ and $t=3$, by Lemma \ref{odd main results}, there is no $3$-SCHGDD of type $(6,1^3)$. When $m=1$, $t\equiv 3,15\ ({\rm mod}\ 18)$ and $t\geq 15$, start from a $3$-SCHGDD of type $(6,(t/3)^3)$, which exists by Lemma \ref{6-SCHGDD-m^3}. Then apply Construction \ref{SCHGDD from SCHGDD} with a $3$-SCHGDD of type $(6,1^{t/3})$ to obtain a $3$-SCHGDD of type $(6,1^t)$.

When $m\equiv 1,5\ ({\rm mod}\ 6)$ and $t\equiv 3,15\ ({\rm mod}\ 18)$, start from a $3$-SCHGDD of type $(6,1^t)$, and apply Construction \ref{SCHGDD from CDM} with a $(3,m)$-CDM to obtain a $3$-SCHGDD of type $(6,m^t)$.

When $t\equiv 0\ ({\rm mod}\ 2)$ and $m\equiv 0\ ({\rm mod}\ 8)$, start from a $3$-SCGDD of type $(m/2)^6$, which exists by Lemma \ref{scgdd}. By Lemma \ref{odd main results}, there exists a $3$-SCHGDD of type $(3,2^t)$ for any even integer $t\geq 4$. Then apply Construction \ref{SCHGDD-recur} to obtain a $3$-SCHGDD of type $(6,m^t)$.

When $t\equiv 0\ ({\rm mod}\ 2)$, $t\neq 8$ and $m=4$, take a strictly cyclic $\{3,5\}$-GDD of type $4^t$ from Lemma \ref{4^t-0,1,2mod3}.  Then apply Construction \ref{SCHGDD-from strictly CGDD} with a $3$-MGDD of type $r^6$ for $r=3,5$, which exists by Theorem \ref{3-HGDD}, to obtain a $3$-SCHGDD of type $(6,4^t)$. When $t=8$ and $m=4$, there is a $3$-SCHGDD of type $(6,4^8)$ by Lemma \ref{6,(2,4)^8}. When $t\equiv 0\ ({\rm mod}\ 2)$ and $m\equiv 4\ ({\rm mod}\ 8)$, start from a $3$-SCHGDD of type $(6,4^t)$, and apply Construction \ref{SCHGDD from CDM} with a $(3,m/4)$-CDM to obtain a $3$-SCHGDD of type $(6,m^t)$.

When $t\equiv 0\ ({\rm mod}\ 4)$ and $m=2$, write $t=2^r u$, where $u\equiv 1\ ({\rm mod}\ 2)$ and $r\geq 2$. Take a $3$-SCHGDD of type $(6,2^{2^r})$ from Lemma \ref{6,2-4^power}. Apply Construction \ref{SCHGDD from CDM} with a $(3,u)$-CDM to obtain a $3$-SCHGDD of type $(6,(2u)^{2^r})$. Then making use of Construction \ref{SCHGDD from SCHGDD} with a $3$-SCHGDD of type $(6,2^u)$, which exists by Lemma \ref{odd main results}, we have a $3$-SCHGDD of type $(6,2^{2^r u})$.

When $t\equiv 2\ ({\rm mod}\ 4)$ and $m=2$, let $I=Z_5\cup\{\infty\}$ and we here construct a $3$-SCHGDD of type $(6,2^t)$ on $I\times Z_{2t}$ with the group set $\{\{i\}\times Z_{2t}:i\in I\}$ and the hole set $\{I\times \{j,t+j\}:0\leq j\leq t-1\}$. Only initial base blocks are listed below, and all other base blocks are obtained by developing these base blocks by $(+1,-)$ modulo $(5,-)$, where $\infty+1=\infty$.

\begin{center}\begin{tabular}{ll}
$\{(0,0),(1,2+i),(2,t+4+2i)\}$, & $i\in[0,t-3]\setminus\{t/2-3,t/2-2\}$;\\
$\{(0,0),(2,3+2i),(\infty,t+2+i)\}$, & $i\in[0,t-4]\setminus\{t/2-2\}$;
\end{tabular}

\begin{tabular}{ll}
$\{(0,0),(1,3t/2-1),(2,t-1)\}$, & $\{(0,0),(2,2t-1),(\infty,1)\}$,\\
$\{(0,0),(1,t/2-1),(\infty,2t-1)\}$, & $\{(0,0),(1,t/2),(\infty,t+1)\}$,\\
$\{(0,0),(1,1),(3,2)\}$, & $\{(0,0),(1,t+1),(3,3)\}$.
\end{tabular}
\end{center}

When $t\equiv 0\ ({\rm mod}\ 2)$ and $m\equiv 2\ ({\rm mod}\ 4)$, start from a $3$-SCHGDD of type $(6,2^t)$, and apply Construction \ref{SCHGDD from CDM} with a $(3,m/2)$-CDM to obtain a $3$-SCHGDD of type $(6,m^t)$. \qed

\begin{Lemma}\label{8-SCHGDD}
If $(t-1)m\equiv 0\ ({\rm mod}\ 6)$ and $t\geq 3$, then there exists a $3$-SCHGDD of type $(8,m^t)$ except possibly when $m\equiv 2,10\ ({\rm mod}\ 12)$ and $t\equiv 7,10\ ({\rm mod}\ 12)$.
\end{Lemma}

\proof When $t$ is odd with the exception of $m\equiv 2,10\ ({\rm mod}\ 12)$ and $t\equiv 7\ ({\rm mod}\ 12)$, the conclusion follows from Lemma \ref{odd main results}. When $m\equiv 6\ ({\rm mod}\ 12)$ and $t\equiv 2\ ({\rm mod}\ 4)$, start from a $4$-SCGDD of type $3^8$, which exists by Lemma $32$ in \cite{wy}. By Lemma \ref{4-SCHGDD}, there exists a $3$-SCHGDD of type $(4,(m/3)^t)$. Then apply Construction \ref{SCHGDD-recur} to obtain a $3$-SCHGDD of type $(8,m^t)$.

When $m\equiv 0\ ({\rm mod}\ 12)$ and $t\equiv 0\ ({\rm mod}\ 2)$, or $m\equiv 6\ ({\rm mod}\ 12)$ and $t\equiv 0\ ({\rm mod}\ 4)$, or $m\equiv 4,8\ ({\rm mod}\ 12)$ and $t\equiv 4\ ({\rm mod}\ 6)$, or $m\equiv 2,10\ ({\rm mod}\ 12)$ and $t\equiv 4\ ({\rm mod}\ 12)$, take a strictly cyclic $3$-GDD of type $m^t$ from Lemma \ref{3-sCGDD}. Then apply Construction \ref{SCHGDD-from strictly CGDD} with a $3$-MGDD of type $3^8$ from Theorem \ref{3-HGDD} to obtain a $3$-SCHGDD of type $(8,m^t)$. \qed

% When $m\equiv 2,10\ ({\rm mod}\ 12)$ and $t\equiv 7,10\ ({\rm mod}\ 12)$, ????????????? \qed

\section{$3$-SCHGDDs of type $(n,m^t)$ for general $n$}

\begin{Lemma}
\label{346pbd} {\rm \cite{abg}}
\begin{enumerate}
\item[$(1)$] There exists a $(v,\{3,4\},1)$-PBD for any integer $v\equiv 0,1\ ({\rm mod}\ 3)$ and $v\geq 3$ with the exception of $v=6$.
\item[$(2)$] There exists a $(v,\{3,4,5\},1)$-PBD for any integer $v\geq 3$ with the exception of $v=6,8$.
%\item[$(3)$] There exists a $(v,\{5,9,13\},1)$-PBD for any integer $v\geq 5$ and $v\equiv 1\ ({\rm mod}\ 4)$ with the exception of $v=17,29,33$.
\end{enumerate}
\end{Lemma}

\begin{Lemma}\label{SCHGDD-n-0-1-no 3}
Let $n\equiv 0,1\ ({\rm mod}\ 3)$ and $n\geq 3$. There exists a $3$-SCHGDD of type $(n,m^t)$ for $(t-1)m\equiv 0\ ({\rm mod}\ 2)$ and $t\geq 3$ with the exception of $(1)$ $n=t=3$ and $m\equiv 0\ ({\rm mod}\ 2)$, $(2)$ $(n,m,t)=(6,1,3)$.
\end{Lemma}

\proof When $n=3,4,6$, the conclusion follows from Lemmas \ref{odd main results}, \ref{4-SCHGDD} and \ref{6-SCHGDD}. When $n\equiv 0,1\ ({\rm mod}\ 3)$ and $n\geq 7$, start from a $\{3,4\}$-SCGDD of type $1^n$, which is also a $(n,\{3,4\},1)$-PBD and exists by Lemma \ref{346pbd}. Take a $3$-SCHGDD of type $(k,m^t)$ for $k\in \{3,4\}$, and then apply Construction \ref{SCHGDD-recur} to obtain a $3$-SCHGDD of type $(n,m^t)$. \qed

To construct $3$-SCHGDDs of type $(n,1^4)$, we introduce the concept of quasi-skew starters. Let $n$ be an odd positive integer and $G$ be an additive abelian group of order $n$. A {\it quasi-skew starter} in $G$ is a set of unordered pairs $\{\{x_i,y_i\}: 1\leq i\leq (n-1)/2\}$, which satisfies the following two properties:
\begin{enumerate}
\item[$(1)$]$\{x_i:1\leq i\leq (n-1)/2\}\cup\{y_i:1\leq i\leq (n-1)/2\}=Z_n\setminus\{0\}$;
\item[$(2)$] $\{\pm(x_i+y_i): 1\leq i\leq (n-1)/2\}=Z_n\setminus\{0\}$.
\end{enumerate}
We remark that there are many kinds of combinatorial configurations named with starter, such as strong starter, skew starter, balanced starter and partitionable starter, etc, which have been used in the construction of various combinatorial designs such as Room squares, Howell designs and Howell rotation. For more details the interested reader may refer to \cite{d}. A quasi-skew starter is said to be a {\em skew starter}, if it further satisfies $\{\pm(x_i-y_i): 1\leq i\leq (n-1)/2\}=Z_n\setminus\{0\}$. That is the reason we use the term quasi-skew starter. The existence of skew starters is far more from solved (cf. \cite{cgz}).

\begin{Lemma}\label{quasi-skew}
There exists a quasi-skew starter in $Z_n$ for all $n\equiv 1\ ({\rm mod}\ 2)$ and $n\geq 7$.
\end{Lemma}

\proof Let $n=2s+1$, where $s\geq 3$. The $s$ pairs of a quasi-skew starter in $Z_{2s+1}$ are listed in Tables $1$ and $2$. \qed

\begin{center}
\begin{table}\small
\begin{tabular}{|c|llllllll|}\hline
$s=3$& $\{1,5\}$ & $\{2,3\}$ & $\{4,6\}$ &&&&& \\\hline
$s=4$& $\{1,5\}$ & $\{3,4\}$ & $\{2,8\}$ & $\{6,7\}$ &&&&\\\hline
$s=5$& $\{1,6\}$ & $\{2,8\}$ & $\{3,5\}$ & $\{4,9\}$ & $\{7,10\}$ &&&\\\hline
$s=6$& $\{1,7\}$ & $\{2,10\}$ & $\{3,8\}$ & $\{4,6\}$ & $\{5,12\}$ & $\{9,11\}$ && \\\hline
$s=7$& $\{1,11\}$ & $\{2,9\}$ & $\{4,5\}$ & $\{6,8\}$ & $\{3,14\}$ & $\{7,13\}$ & $\{10,12\}$&\\\hline
$s=8$& $\{1,14\}$ & $\{2,10\}$ & $\{3,11\}$ & $\{4,6\}$ & $\{7,9\}$ & $\{5,16\}$ & $\{8,15\}$ & $\{12,13\}$\\\hline
$s=10$& $\{1,15\}$ & $\{2,16\}$ & $\{3,14\}$ & $\{4,9\}$ & $\{6,8\}$ & $\{7,13\}$ & $\{5,18\}$ & $\{10,17\}$\\ & $\{11,19\}$ & $\{12,20\}$&&&&&& \\\hline
$s=11$ & $\{1,13\}$ & $\{2,18\}$ & $\{4,14\}$ & $\{5,11\}$ & $\{6,16\}$ & $\{7,12\}$ & $\{8,9\}$ & $\{3,22\}$\\ &$\{10,21\}$ & $\{15,20\}$ & $\{17,19\}$&&&&& \\\hline
$s=13$& $\{1,15\}$ & $\{2,18\}$ & $\{4,13\}$ & $\{5,19\}$ & $\{7,11\}$ & $\{8,14\}$ & $\{10,16\}$ & $\{3,26\}$\\ &$\{6,25\}$ & $\{9,24\}$ & $\{12,23\}$ & $\{17,22\}$ & $\{20,21\}$&&& \\\hline
$s=14$& $\{1,15\}$ & $\{2,18\}$ & $\{4,20\}$ & $\{5,13\}$ & $\{6,22\}$ & $\{7,19\}$ & $\{9,16\}$ & $\{10,12\}$\\ &$\{3,28\}$ & $\{8,27\}$ & $\{11,26\}$ & $\{14,25\}$ & $\{17,24\}$ & $\{21,23\}$&&\\\hline
\end{tabular}
\caption{quasi-skew starters in $Z_{2s+1}$ for small $s$}
\end{table}
\end{center}

\begin{center}
\begin{table}\small
\begin{tabular}{|c|l|}\hline
& $\{i,i+s\}$, $i\in[1,s/2-1]\setminus\{s/8+1,s/4\}$;\\
$s\equiv 0\ ({\rm mod}\ 8)$ & $\{i+s/2+1,i-s/2-1\}$, $i\in[2,s/2-1]\setminus\{3s/8+1\}$;\\
$s\geq 16$ &$\{s/2,s/2+1\},$  $\{5s/4,3s/2\},$  $\{s/8+1,15s/8+1\},$\\
&$\{3s/2+1,2s\},$ $\{s/4,s/2+2\},$  $\{7s/8+2,9s/8+1\}.$\\\hline

& $\{i,i+s\}$, $i\in[1,(s-3)/2]\setminus\{(s-1)/8,(s-1)/4\}$;\\
$s\equiv 1\ \{{\rm mod}\ 8\}$ & $\{i+(s+3)/2,i-(s+3)/2+1\}$, $i\in[1,(s-3)/2]\setminus\{(3s-11)/8\}$;\\
$s\geq9$ & $\{(s-1)/8,(15s-7)/8\},$ $\{(s+3)/2,(s-1)/4\},$ $\{(s-1)/2,(s+1)/2\},$\\
& $\{(7s+1)/8,(9s-1)/8\},$ $\{(3s-1)/2,(5s-1)/4\},$ $\{(3s+1)/2,2s\}.$\\\hline

& $\{i,i+s\}$, $i\in[2,s/2-1]\setminus\{(s+6)/8\}$;\\
$s\equiv 2\ ({\rm mod}\ 8\}$ & $\{i+s/2+1,i-s/2-1\}$, $i\in[2,s/2-1]\setminus\{(s-2)/4,(3s+10)/8\}$;\\
$s\geq18$ & $\{s/2,s/2+2\},$  $\{(7s+18)/8,(9s+6)/8\},$  $\{(s+6)/8,(15s+10)/8\},$ \\
& $\{s/2+1,(3s+2)/4\},$ $\{3s/2,(7s-2)/4\},$  $\{1,3s/2+1\},$  $\{s+1,2s\}.$  \\\hline

& $\{i,i+s\}$, $i\in[1,(s-3)/2]\setminus\{(s-3)/8\}$;\\
$s\equiv 3\ ({\rm mod}\ 8)$ & $\{i+(s+3)/2,i-(s+3)/2+1\}$, $i\in[1,(s-3)/2]\setminus\{(3s-9)/8,(s+1)/4\}$;\\
$s\geq19$ & $\{(s-3)/8,(15s-5)/8\},$  $\{(s+1)/2,(3s+7)/4\},$  $\{(3s+1)/2,2s\},$\\
& $\{(7s+3)/8,(9s-3)/8\},$  $\{(3s-1)/2,(7s+3)/4\},$  $\{(s-1)/2,(s+3)/2\}.$\\\hline

& $\{i,i+s\}$, $i\in[1,s/2-1]\setminus\{(s+4)/8\}$;\\
$s\equiv 4\ ({\rm mod}\ 8)$ & $\{i+s/2+1,i-s/2-1\}$, $i\in[2,s/2-1]\setminus\{(3s-12)/8,s/4+1\}$;\\
$s\geq12$ & $\{(s+4)/8,(15s-12)/8\},$ $\{3s/2+1,7s/4+1\},$ $\{3s/2,2s\},$\\
& $\{(7s-4)/8,(9s+4)/8\},$ $\{s/2+2,3s/4+2\},$ $\{s/2,s/2+1\}.$\\\hline

& $\{i,i+s\}$, $i\in[1,(s-3)/2]\setminus\{(s+11)/8,(s-1)/4\}$;\\
$s\equiv 5\ ({\rm mod}\ 8)$ & $\{i+(s+3)/2,i-(s+1)/2\}$, $i\in[1,(s-3)/2]\setminus\{(3s+1)/8\}$;\\
$s\geq 21$ & $\{(3s-1)/2,(5s-1)/4\},$ $\{(9s+11)/8,(7s+13)/8\},$ $\{(s-1)/2,(s+1)/2\},$\\
& $\{(s+3)/2,(s-1)/4\},$ $\{(3s+1)/2,2s\},$ $\{(s+11)/8,(15s+5)/8\}.$\\\hline

& $\{i,i+s\}$, $i\in[2,s/2-1]\setminus\{(s-6)/8\}$;\\
$s\equiv 6\ ({\rm mod}\ 8)$ & $\{i+s/2+1,i-s/2-1\}$, $i\in[2,s/2-1]\setminus\{(3s-2)/8,(s-2)/4\}$;\\
$s\geq 22$ & $\{(s-6)/8,(15s-2)/8\},$ $\{3s/2,(7s-2)/4\},$ $\{1,3s/2+1\},$ $\{s+1,2s\},$\\
& $\{(7s+6)/8,(9s-6)/8\},$ $\{s/2+1,(3s+2)/4\},$ $\{s/2,s/2+2\}.$\\\hline

& $\{i,i+s\}$, $i\in[1,(s-3)/2]\setminus\{(s+9)/8\}$;\\
$s\equiv 7\ ({\rm mod}\ 8)$ & $\{i+(s+3)/2,i+(3s+1)/2\}$, $i\in[1,(s-3)/2]\setminus\{(s+1)/4,(3s+3)/8\}$;\\
$s\geq15$ & $\{(9s+9)/8,(7s+15)/8\},$ $\{(3s+1)/2,2s\},$ $\{(3s-1)/2,(7s+3)/4\},$\\
& $\{(s-1)/2,(s+3)/2\},$ $\{(s+9)/8,(15s+7)/8\},$ $\{(s+1)/2,(3s+7)/4\}.$\\\hline
\end{tabular}
\caption{quasi-skew starters in $Z_{2s+1}$ for general $s$}
\end{table}
\end{center}

\begin{Lemma}\label{n-1,1^4} Let $n\equiv 1\ ({\rm mod}\ 2)$ and $n\geq 7$. There exists a $3$-SCHGDD of type $(n,1^4)$.
\end{Lemma}

\proof We here construct a $3$-SCHGDD of type $(n,1^4)$ on $Z_n\times Z_4$ with the group set $\{\{i\}\times Z_{4}:i\in Z_n\}$ and the hole set $\{Z_n\times \{j\}: j\in Z_4\}$.
The required $n(n-1)/2$ base blocks are
\begin{center}
$\{(i,0),(x_r+i,1),(x_r+y_r+i,2)\}$, $1\leq r\leq (n-1)/2$, $0\leq i\leq n-1$,
\end{center}
where $\{\{x_r,y_r\}: 1\leq r\leq (n-1)/2\}$ is a quasi-skew starter in $Z_n$. \qed

\begin{Lemma}\label{n-1-3,m^t-0mod4} Let $n\equiv 1,3\ ({\rm mod}\ 6)$, $n\geq 7$ and $m\equiv 1\ ({\rm mod}\ 2)$. There exists a $3$-SCHGDD of type $(n,m^t)$ for any $t\equiv 0\ ({\rm mod}\ 4)$ except possibly when $t=8$.
\end{Lemma}

\proof When $m=1$ and $t=4$, the conclusion follows from Lemma \ref{n-1,1^4}. When $m=1$, $t\equiv 0\ ({\rm mod}\ 4)$ and $t\geq 12$, write $t=4t_1$. Take a $3$-SCHGDD of type $(n,4^{t_1})$ from Lemma \ref{SCHGDD-n-0-1-no 3}. Apply Construction \ref{SCHGDD from SCHGDD} with a $3$-SCHGDD of type $(n,1^4)$ to obtain a $3$-SCHGDD of type $(n,1^{4t_1})$. When $m\equiv 1\ ({\rm mod}\ 2)$ and $t\equiv 0\ ({\rm mod}\ 4)$, start from a $3$-SCHGDD of type $(n,1^t)$, and apply Construction \ref{SCHGDD from CDM} with a $(3,m)$-CDM to obtain a $3$-SCHGDD of type $(n,m^t)$. \qed

\begin{Lemma}{\rm\cite{wy}}\label{SCGDD from SCHGDD}
Suppose that a $k$-SCHGDD of type $(n,m^t)$ and a $k$-SCGDD of type $m^n$ exist. Then a $k$-SCGDD of type $(mt)^n$ exists.
\end{Lemma}

\begin{Lemma}\label{n-3-7-12-no} Let $n\equiv 3,7\ ({\rm mod}\ 12)$. There is no $3$-SCHGDD of type $(n,m^t)$ for any positive integer $m\equiv 1\ ({\rm mod}\ 2)$ and $t\equiv 2\ ({\rm mod}\ 4)$.
\end{Lemma}

\proof Suppose that there was a $3$-SCHGDD of type $(n,m^t)$ for $m\equiv 1\ ({\rm mod}\ 2)$ and $t\equiv 2\ ({\rm mod}\ 4)$. Applying Construction \ref{SCHGDD from SCHGDD} with a $3$-SCHGDD of type $(n,1^m)$, which exists by Lemma \ref{SCHGDD-n-0-1-no 3}, we would have a $3$-SCHGDD of type $(n,1^{mt})$. However, when $n\equiv 3,7\ ({\rm mod}\ 12)$, Lemma \ref{SCGDD from SCHGDD} shows that a $3$-SCHGDD of type $(n,1^{mt})$ implies existence of a $3$-SCGDD of type $(mt)^n$. Since $mt\equiv 2\ ({\rm mod}\ 4)$, this SCGDD cannot exist by Lemma \ref{scgdd}, a contradiction. \qed

\begin{Lemma}\label{SCHGDD-n-2-(1)}
Let $n\equiv 2\ ({\rm mod}\ 3)$ and $n\geq 5$. There exists a $3$-SCHGDD of type $(n,m^t)$ for $(t-1)m\equiv 0\ ({\rm mod}\ 6)$ and $t\geq 3$ except possibly when $n=8$, $m\equiv 2,10\ ({\rm mod}\ 12)$ and $t\equiv 7,10\ ({\rm mod}\ 12)$.
\end{Lemma}

\proof When $t$ is odd, the conclusion follows from Lemma \ref{odd main results}. When $n=5,8$, the conclusion follows from Lemmas \ref{5-SCHGDD} and \ref{8-SCHGDD}. Assume that $n\equiv 2\ ({\rm mod}\ 3)$, $n\geq 11$ and $t$ is even. Start from a $\{3,4,5\}$-SCGDD of type $1^n$, which is also a $(n,\{3,4,5\},1)$-PBD and exists by Lemma \ref{346pbd}. For $(t-1)m\equiv 0\ ({\rm mod}\ 6)$ and $t\geq 4$, take a $3$-SCHGDD of type $(k,m^t)$ for $k\in \{3,4,5\}$ from Lemmas \ref{odd main results}, \ref{4-SCHGDD} and \ref{5-SCHGDD}. Then apply Construction \ref{SCHGDD-recur} to obtain a $3$-SCHGDD of type $(n,m^t)$. \qed

\begin{Lemma}\label{n-5,m^t-4mod12} Let $n\equiv 5\ ({\rm mod}\ 6)$ and $n\geq 5$. There exists a $3$-SCHGDD of type $(n,m^t)$ for any positive integer $m\equiv 1\ ({\rm mod}\ 2)$ and $t\equiv 4\ ({\rm mod}\ 12)$ with the exception of $(n,m,t)=(5,1,4)$.
\end{Lemma}

\proof When $n=5$, the conclusion follows from Lemma \ref{5-SCHGDD}. When $n\equiv 5\ ({\rm mod}\ 6)$ and $n\geq 11$, if $m=1$ and $t=4$, the conclusion follows from Lemma \ref{n-1,1^4}. If $m=1$, $t\equiv 4\ ({\rm mod}\ 12)$ and $t\geq 16$, write $t=4t_1$. Take a $3$-SCHGDD of type $(n,4^{t_1})$ from Lemma \ref{SCHGDD-n-2-(1)}. Apply Construction \ref{SCHGDD from SCHGDD} with a $3$-SCHGDD of type $(n,1^4)$ to obtain a $3$-SCHGDD of type $(n,1^{4t_1})$. If $m\equiv 1\ ({\rm mod}\ 2)$ and $t\equiv 4\ ({\rm mod}\ 12)$, start from a $3$-SCHGDD of type $(n,1^t)$, and apply Construction \ref{SCHGDD from CDM} with a $(3,m)$-CDM to obtain a $3$-SCHGDD of type $(n,m^t)$. \qed

\begin{Lemma}\label{n-5,3^4} Let $n\equiv 5\ ({\rm mod}\ 6)$. There exists a $3$-SCHGDD of type $(n,3^4)$.
\end{Lemma}

\proof We here construct a $3$-SCHGDD of type $(n,3^4)$ on $Z_n\times Z_{12}$ with the group set $\{\{i\}\times Z_{12}:i\in Z_n\}$ and the hole set $\{Z_n\times \{j,4+j,8+j\}:0\leq j\leq 3\}$.
All the $3n(n-1)/2$ base blocks can be obtained from the following $3(n-1)/2$ initial blocks by $(+1\ {\rm mod}\ n,-)$. When $n=5$, the $6$ initial blocks are
\begin{center}
\begin{tabular}{lll}
$\{(0,0),(1,1),(4,3)\}$, & $\{(0,0),(2,1),(3,3)\}$, &
$\{(0,0),(2,2),(3,5)\}$, \\ $\{(0,0),(4,2),(1,5)\}$, &
$\{(0,0),(3,1),(2,6)\}$, & $\{(0,0),(4,1),(1,6)\}$. \\
\end{tabular}
\end{center}
When $n=11$, the $15$ initial blocks are
\begin{center}
\begin{tabular}{lll}
$\{(0,0),(6,1),(1,3)\}$, & $\{(0,0),(1,2),(3,5)\}$, & $\{(0,0),(1,1),(3,6)\}$,\\
$\{(0,0),(7,1),(3,3)\}$, & $\{(0,0),(5,2),(4,5)\}$, & $\{(0,0),(4,1),(5,6)\}$,\\
$\{(0,0),(8,1),(5,3)\}$, & $\{(0,0),(4,2),(8,5)\}$, & $\{(0,0),(5,1),(1,6)\}$,\\
$\{(0,0),(9,1),(7,3)\}$, & $\{(0,0),(3,2),(9,5)\}$, & $\{(0,0),(3,1),(9,6)\}$,\\
$\{(0,0),(10,1),(9,3)\}$, & $\{(0,0),(2,2),(10,5)\}$, & $\{(0,0),(2,1),(7,6)\}$.\\
\end{tabular}
\end{center}
When $n\equiv 5\ ({\rm mod}\ 6)$ and $n\geq 17$, write $n=6s+5$ and $s\geq 2$. The $9s+6$ initial blocks are divided into two parts. The first part consists of the following $6s+7$ initial blocks:
\begin{center}
\begin{tabular}{ll}
$\{(0,0),(1+i,1),(2+2i,3)\}$,& $i\in[0,3s+1]$;\\
$\{(0,0),(6s+4-i,2),(2+i,5)\}$,& $i\in[0,3s]$;\\
\end{tabular}
\begin{tabular}{ll}
$\{(0,0),(3s+3,2),(3s+4,5)\}$, & $\{(0,0),(3s+3,1),(1,6)\}$,\\
$\{(0,0),(6s+2,1),(6s+3,6)\}$, & $\{(0,0),(6s+4,1),(6s+2,6)\}$.\\
\end{tabular}
\end{center}
The second part consists of the following $3s-1$ initial blocks:
\begin{itemize}
\item if $s\equiv 0\ ({\rm mod}\ 2)$ and $s\geq 2$, then
\begin{center}
\begin{tabular}{lll}
$\{(0,0),(3s+4+i,1),(4+2i,6)\}$, &  $i\in[0,3s/2-3]$;\\
$\{(0,0),(9s/2+5+3i,1),(3s+8+6i,6)\}$, & $i\in[0,s/2-2]$, & (null if $s=2$);\\
$\{(0,0),(9s/2+7+3i,1),(3s+12+6i,6)\}$, & $i\in[0,s/2-2]$, & (null if $s=2$);\\
$\{(0,0),(9s/2+9+3i,1),(3s+10+6i,6)\}$, & $i\in[0,s/2-2]$, & (null if $s=2$);\\
\end{tabular}
\begin{tabular}{ll}
$\{(0,0),(9s/2+2,1),(3s+2,6)\}$, & $\{(0,0),(9s/2+3,1),(3s+1,6)\}$,\\
$\{(0,0),(9s/2+4,1),(3s+6,6)\}$, & $\{(0,0),(9s/2+6,1),(3s+5,6)\}$;\\
\end{tabular}
\end{center}

\item if $s\equiv 1\ ({\rm mod}\ 2)$ and $s\geq 3$, then
\begin{center}
\begin{tabular}{lll}
$\{(0,0),(3s+4+i,1),(4+2i,6)\}$, &  $i\in[0,(3s-5)/2]$;\\
$\{(0,0),((9s+11)/2+3i,1),(3s+9+6i,6)\}$, & $i\in[0,(s-3)/2]$;\\
$\{(0,0),((9s+15)/2+3i,1),(3s+7+6i,6)\}$, & $i\in[0,(s-3)/2]$;\\
$\{(0,0),((9s+13)/2+3i,1),(3s+11+6i,6)\}$, & $i\in[0,(s-5)/2]$, & (null if $s=3$);\\
\end{tabular}
\begin{tabular}{ll}
$\{(0,0),((9s+5)/2,1),(3s+4,6)\}$, & $\{(0,0),((9s+7)/2,1),(3s+2,6)\}$,\\
$\{(0,0),((9s+9)/2,1),(3s+5,6)\}$.\\
\end{tabular}
\end{center}
\end{itemize}
%\qed

%\begin{Lemma}\label{n-5,3^8} Let $n\equiv 5\ ({\rm mod}\ 6)$. There exists a $3$-SCHGDD of type $(n,3^8)$.
%\end{Lemma}

%\proof ???????????????????????

\begin{Lemma}\label{n-5,m^t-0mod4} Let $n\equiv 5\ ({\rm mod}\ 6)$ and $m\equiv 3\ ({\rm mod}\ 6)$. There exists a $3$-SCHGDD of type $(n,m^t)$ for any $t\equiv 0\ ({\rm mod}\ 4)$ except possibly when $t=8$.
\end{Lemma}

\proof When $m=3$ and $t=4$, the conclusion follows from Lemma \ref{n-5,3^4}. When $m=3$, $t\equiv 0\ ({\rm mod}\ 4)$ and $t\geq 12$, write $t=4t_1$. Take a $3$-SCHGDD of type $(n,12^{t_1})$ from Lemma \ref{SCHGDD-n-2-(1)}. Apply Construction \ref{SCHGDD from SCHGDD} with a $3$-SCHGDD of type $(n,3^4)$ to obtain a $3$-SCHGDD of type $(n,3^{4t_1})$. When $m\equiv 3\ ({\rm mod}\ 6)$ and $t\equiv 0\ ({\rm mod}\ 4)$, start from a $3$-SCHGDD of type $(n,3^t)$, and apply Construction \ref{SCHGDD from CDM} with a $(3,m/3)$-CDM to obtain a $3$-SCHGDD of type $(n,m^t)$. \qed

%\begin{Lemma}\label{n-5,m^t-2mod4} Let $n\equiv 5\ ({\rm mod}\ 6)$. There exists a $3$-SCHGDD of type $(n,m^t)$ for any positive integer $m\equiv 3\ ({\rm mod}\ 6)$ and $t\equiv 2\ ({\rm mod}\ 4)$.
%\end{Lemma}

%\proof ?????????????????????????? \qed

%\begin{Lemma}\label{n-5,m^t-10mod12} Let $n\equiv 5\ ({\rm mod}\ 6)$. There exists a $3$-SCHGDD of type $(n,m^t)$ for any positive integer $m\equiv 1,5\ ({\rm mod}\ 6)$ and $t\equiv 10\ ({\rm mod}\ 12)$.
%\end{Lemma}

%\proof ?????????????????????????? \qed

Now combining the results of Lemmas \ref{odd main results}, \ref{nece}, \ref{SCHGDD-n-0-1-no 3}, \ref{n-1-3,m^t-0mod4}, \ref{n-3-7-12-no}-\ref{n-5,m^t-4mod12} and \ref{n-5,m^t-0mod4}, we obtain the main theorem in this paper as follows.

\begin{Theorem}\label{main theorem}
There exists a $3$-SCHGDD of type $(n,m^t)$ if and only if $n,t\geq 3$, $(t-1)(n-1)m\equiv 0\ ({\rm mod}\ 2)$ and $(t-1)n(n-1)m\equiv 0\ ({\rm mod}\ 6)$ except when
\begin{enumerate}
\item[$(1)$] $n\equiv 3,7\ ({\rm mod}\ 12)$, $m\equiv 1\ ({\rm mod}\ 2)$ and $t\equiv 2\ ({\rm mod}\ 4)$;
\item[$(2)$] $n=3$, $m\equiv 1\ ({\rm mod}\ 2)$ and $t\equiv 0\ ({\rm mod}\ 2)$;
\item[$(3)$] $n=t=3$ and $m\equiv 0\ ({\rm mod}\ 2)$;
\item[$(4)$] $(n,m,t)\in\{(5,1,4),(6,1,3)\}$;
\end{enumerate}
and possibly when
\begin{enumerate}
\item[$(1)$] $n=8$, $m\equiv 2,10\ ({\rm mod}\ 12)$ and $t\equiv 7,10\ ({\rm mod}\ 12)$;
\item[$(2)$] $t=8$, either $m\equiv 1\ ({\rm mod}\ 2)$, $n\equiv 1,3\ ({\rm mod}\ 6)$ and $n\geq 7$, or $m\equiv 3\ ({\rm mod}\ 6)$ and $n\equiv 5\ ({\rm mod}\ 6)$;
\item[$(3)$] $n\equiv 1,9\ ({\rm mod}\ 12)$, $m\equiv 1\ ({\rm mod}\ 2)$ and $t\equiv 2\ ({\rm mod}\ 4)$;
\item[$(4)$] $n\equiv 5\ ({\rm mod}\ 6)$, either $m\equiv 3\ ({\rm mod}\ 6)$ and $t\equiv 2\ ({\rm mod}\ 4)$, or $m\equiv 1,5\ ({\rm mod}\ 6)$ and $t\equiv 10\ ({\rm mod}\ 12)$.
\end{enumerate}
\end{Theorem}

\section{Applications}

\subsection{Two-dimensional balanced sampling plans}

In environmental and ecological populations, neighboring units within a finite population, spatially or sequentially ordered, may provide similar information. In an attempt to obtain the most informative picture of a population, one desires a sample avoiding the selection of adjacent units. Balanced sampling plans excluding adjacent units have been proposed as a means of achieving such a goal (see for example \cite{hrs,stu}).

There are many types of two-dimensional balanced sampling plans excluding adjacent units, which have different adjacency scheme \cite{wr}. Here two-dimensional means the set of $nm$ units (called {\em points}) within populations, say $Z_n\times Z_m$, is arranged in two dimensions naturally.

Given $(x,y)\in Z_n\times Z_m$, the points $(i,y)$ and $(x,j)$ for any $i\in Z_n$ and $j\in Z_m$ are said to be {\em row-column-mates} of the point $(x, y)$. A {\em two-dimensional balanced sampling plan excluding row-column-mates} is a pair $(X,{\cal B})$, where $X=Z_n\times Z_m$ and $\cal B$ is a collection of $k$-subsets of $X$ (called {\em blocks}) such that any two points that are row-column-mates do not appear in any block while any two points that are not row-column-mates appear in exactly one block. Obviously it is just a $k$-HGDD of type $(n,1^m)$, or equivalently, a $k$-MGDD of type $m^n$.

So our results for semi-cyclic HGDDs can be translated simply using the language of two-dimensional balanced sampling plans excluding row-column-mates. It is meaningful to do this because the resulting balanced sampling plans admit good algebraic structures and have clear advantage over those with no algebraic structures in the identification of the supports.

Another kind of balanced sampling plans we shall mention are balanced sampling plan excluding contiguous units. For $(x,y)\in Z_n\times Z_m$, the points
$(x-1,y)$, $(x+1,y)$, $(x,y-1)$, and $(x,y+1)$ (reducing the sums modulo $n$ and $m$ in the first and second coordinates, respectively) are said to be {\em $2$-contiguous} to the point $(x, y)$. A {\em two-dimensional balanced sampling plan excluding contiguous units} is a pair $(X,{\cal B})$, where $X=Z_n\times Z_m$ and $\cal B$ is a collection of $k$-subsets of $X$ (called {\em blocks}) such that any two $2$-contiguous points do not appear in any block while any two points that are not $2$-contiguous appear in exactly one block. Such a design is denoted by a $2$-BSEC$(n,m,k,1)$ (cf. \cite{bcrw}).

When $n=1$ or $m=1$, a $2$-BSEC can be seen as a balanced sampling plan excluding contiguous units with only one dimension, and is often said to be {\em one-dimensional}, simply written as $1$-BSEC$(m,k,1)$ or $1$-BSEC$(n,k,1)$. In this case each point has only two $2$-contiguous points.

\begin{Theorem}\label{3-BSEC}
\begin{enumerate}
\item[$(1)$]{\rm \cite{cl}} There exists a $1$-BSEC$(m,3,1)$ if and only if $m\geq 9$ and $m\equiv 3\ ({\rm mod}\ 6)$.
\item[$(2)$]{\rm \cite{bcrw}} Let $n,m\geq 3$. There exists a $2$-BSEC$(n,m,3,1)$ if and only if $n$ and $m$ are odd and either $n\equiv m\equiv 3\ ({\rm mod}\ 6)$ or $n\not\equiv m\ ({\rm mod}\ 6)$.
\end{enumerate}
\end{Theorem}

An {\em automorphism} of a $2$-BSEC $(X,{\cal B})$ is a permutation on $X$ leaving ${\cal B}$ invariant. A $2$-BSEC$(n,m,k,1)$ is said to be {\em semi-cyclic} if it admits the automorphism $\pi:(x,y)\longmapsto(x,y+1)$ (mod $(-,m)$), $(x,y)\in Z_n\times Z_m$. Let $H$ be the cyclic group generated by $\pi$ under the compositions of permutations. Then all blocks of the $2$-BSEC can be partitioned into some block orbits under $H$. Choose any fixed block from each block orbit and then call it a {\em base block} of this $2$-BSEC under $H$. Note that the automorphism can be also taken as $(x,y)\longmapsto(x+1,y)$ (mod $(n,-)$), but by the symmetry of $n$ and $m$, there is no difference between them essentially. So we always use the former one.

When $n=1$, a semi-cyclic $2$-BSEC$(1,m,k,1)$ is often said to be a {\em cyclic} $1$-BSEC$(m,k,1)$. Wei \cite{wei1} established the necessary and sufficient conditions for the existence of cyclic $1$-BSECs for block size three.

\begin{Theorem}\label{3-CBSEC}{\rm \cite{wei1}}
There exists a cyclic $1$-BSEC$(m,3,1)$ if and only if $m\geq 9$ and $m\equiv 3\ ({\rm mod}\ 6)$.
\end{Theorem}

As an application of $3$-SCHGDDs, we here present an infinite family of semi-cyclic $2$-BSEC$(n,m,3,1)$s.

\begin{Theorem}\label{n-3;m-3}
There exists a semi-cyclic $2$-BSEC$(n,m,3,1)$ for any $m,n\geq 9$ and $m,n\equiv 3\ ({\rm mod}\ 6)$.
\end{Theorem}

\proof By Theorem \ref{main theorem}, there exists a $3$-SCHGDD of type $(n,1^m)$ $(X,{\cal G},{\cal H},{\cal A})$ for $m,n\geq 9$ and $m,n\equiv 3\ ({\rm mod}\ 6)$, where $X=Z_n\times Z_m$, ${\cal G}=\{G_i:i\in Z_n\}=\{\{i\}\times Z_m:i\in Z_n\}$ and ${\cal H}=\{H_j:j\in Z_m\}=\{Z_n\times \{j\}:j\in Z_m\}$. Denote the set of its base blocks by ${\cal A}^*$. For each group $G_i$ of the $3$-SCHGDD, construct a cyclic $1$-BSEC$(m,3,1)$ $(G_i,{\cal B}_i)$, which exists by Theorem \ref{3-CBSEC}. Denote the set of its base blocks by ${\cal B}^*_i$. For the hole $H_0$ of the $3$-SCHGDD, construct a $1$-BSEC$(n,3,1)$ $(H_0,{\cal C})$, which exists by Theorem \ref{3-BSEC}. Then it is readily checked that ${\cal A}^*\cup(\bigcup_{i=1}^n {\cal B}^*_i)\cup{\cal C}$ is the set of base blocks of the desired semi-cyclic $2$-BSEC$(n,m,3,1)$.  \qed

By counting the necessary conditions for the existence of semi-cyclic $2$-BSEC$(n,m,3,1)$s, it seems that a complete solution for this problem will rely heavily on direct constructions. So finding more recursive constructions to obtain semi-cyclic $2$-BSECs (not only focus on $k=3$) is an interesting problem.

\subsection{Two-dimensional optical orthogonal codes}

An optical orthogonal code is a family of sequences with good auto-
and cross-correlation properties. Its study has been motivated by an
application in an optical code-division multiple access (OCDMA)
system (see for example \cite{csw,dx,fm,sb}).

Let $n$, $m$, $k$ and $\lambda$ be positive
integers. A {\em two-dimensional $(n\times m,k,\lambda)$ optical
orthogonal code} (briefly $2$-D $(n\times m,k,\lambda)$-OOC),
$\cal{C}$, is a family of $n\times m$ $(0,1)$-matrices (called {\em
codewords}) of Hamming weight $k$ satisfying the following two properties
(the arithmetic $j+r$ is reduced modulo $m$):

\begin{enumerate}
\item[(1)] The Auto-Correlation Property: for each matrix
${\mathbf{A}}=(a_{ij})_{n\times m}\in\cal{C}$ and each
integer $r$, $r\not\equiv 0\ ({\rm mod}\ m)$,
$\sum_{i=0}^{n-1}\sum_{j=0}^{m-1}a_{ij}a_{i,j+r}\leq\lambda;$

\item[(2)] The Cross-Correlation Property: for each matrix
${\mathbf{A}}=(a_{ij})_{n\times m}\in\cal{C}$, ${\mathbf{B}}=(b_{ij})_{n\times m}\in\cal{C}$ with ${\mathbf{A}}\neq {\mathbf{B}}$, and each integer $r$,
$\sum_{i=0}^{n-1}\sum_{j=0}^{m-1}a_{ij}b_{i,j+r}\leq\lambda.$
\end{enumerate}

A convenient way of viewing optical orthogonal codes is from a set-theoretic perspective. Consider a $2$-D $(n\times m,k,\lambda)$-OOC, $\cal C$. For each $n\times m$
$(0,1)$-matrix $M\in\cal C$, whose rows are indexed by $I_n$ and
columns are indexed by $Z_m$. Construct a $k$-subset $B_M$ of $I_n\times Z_m$ such that $(i,j)\in B_M$ if and only if $M$'s $(i,j)$ cell equals $1$. Then $\{B_M: M\in {\cal C}\}$ is a set-theoretic representation of the $2$-D $(n\times m,k,\lambda)$-OOC.

\begin{Example} \label{example-2-4} The following two $(0,1)$-matrices constitute a $2$-D $(2\times 4,3,1)$-OOC.

\begin{center}
\begin{tabular}{ll}
$\left(
  \begin{array}{cccc}
    1 & 1&  0& 0 \\
    1 & 0&  0& 0 \\
  \end{array}
\right),$&
$\left(
  \begin{array}{cccc}
    0 & 0& 0& 1 \\
    1 & 1& 0& 0 \\
  \end{array}
\right).$
\end{tabular}
\end{center}

\noindent In set notation, it consists of two $3$-subsets of $I_2\times Z_4$:
$\{(0,0),(0,1),(1,0)\}$ and $\{(1,0),(1,1)$, $(0,3)\}$.

\end{Example}

The number of codewords of a $2$-D OOC is called the {\em size} of the $2$-D OOC. From a practical point of view, a code with a large size is required. For fixed values of $n$, $m$, $k$ and $\lambda$, the largest possible size of a $2$-D
$(n\times m,k,\lambda)$-OOC is denoted by $\Phi(n\times m,k,\lambda)$. A $2$-D $(n\times m,k,\lambda)$-OOC with $\Phi(n\times m,k,\lambda)$ codewords is said to be {\em
optimal}. Generally speaking, it is difficult to determine the exact value of $\Phi(n\times m,k,\lambda)$. Based on the Johnson bound for constant weight
codes, the size of a $2$-D $(u\times v,k,\lambda)$-OOC is upper bounded
\cite{yk} by
$
 \Phi(u\times v,k,\lambda)\leq J(u\times v,k,\lambda),
$
where
$J(u\times v,k,\lambda)=\lfloor\frac{u}{k}\lfloor\frac{uv-1}{k-1}
\lfloor\frac{uv-2}{k-2}\lfloor\cdots\lfloor\frac{uv-\lambda}{k-\lambda}
\rfloor\cdots\rfloor\rfloor\rfloor\rfloor.$ So the $2$-D $(2\times 4,3,1)$-OOC shown in Example \ref{example-2-4} is optimal.

\begin{Lemma}
\label{relation} {\rm \cite{am1}} Suppose that there is a $2$-D $(n\times m,k,\lambda)$-OOC with $b$ codewords. Then for any integer factorization
$m=m_1m_2$, there is a $2$-D $(nm_1\times m_2,k,\lambda)$-OOC with $m_1b$ codewords.
\end{Lemma}

\begin{Lemma}\label{n-4} There exists an optimal $2$-D $(n\times 4,3,1)$-OOC with $J(n\times 4,3,1)=n(2n-1)/3$ codewords for any positive integer $n\equiv 0,2\ ({\rm mod}\ 6)$.
\end{Lemma}

\proof When $n=2$, the conclusion follows from Example \ref{example-2-4}. When $n=6$, a $2$-D $(3\times 8,3,1)$-OOC with $11$ codewords is explicit constructed in the proof of Lemma $5.9$ in \cite{wcf}. Then apply Lemma \ref{relation} to obtain the required $2$-D $(6\times 4,3,1)$-OOC with $22$ codewords.

Assume that $n\equiv 0,2\ ({\rm mod}\ 6)$ and $n\geq 8$. By Lemma \ref{scgdd}, there exists a $3$-SCGDD of type $4^{n/2}$, which is defined on $I_{n/2}\times Z_4$ and with base block set ${\cal B}^*$. It is well known that a $3$-GDD of type $2^3$ exists. For each $B=\{(i_x,x),(i_y,y),(i_z,z)\}\in{\cal B}^*$, construct a $3$-GDD of type $2^3$ on the point set $B\cup \{(i_x+n/2,x),(i_y+n/2,y),(i_z+n/2,z)\}$ with the group set $\{\{(i_r,r),(i_r+n/2,r)\}:r\in\{x,y,z\}\}$. Denote the set of its blocks by ${\cal A}_B$. For each $j\in I_{n/2}$, construct an optimal $2$-D $(2\times 4,3,1)$-OOC with $2$ codewords, ${\cal C}_j$, on the set $\{j,j+n/2\}\times Z_4$. Now we construct the required $2$-D $(n\times 4,3,1)$-OOC on $I_n\times Z_4$. It is readily checked that $(\bigcup_{B\in {\cal B}^*}{\cal A}_B)\cup(\bigcup_{j\in I_{n/2}}{\cal C}_j)$ is the set of the required $n(2n-1)/3$ codewords. \qed

\begin{Theorem}\label{OOC-1} There exists an optimal $2$-D $(n\times m,3,1)$-OOC with $J(n\times m,3,1)$ codewords for any $n\equiv 0,2\ ({\rm mod}\ 6)$ and $m\equiv 4\ ({\rm mod}\ 12)$.
\end{Theorem}

\proof When $m=4$, the conclusion follows from Lemma \ref{n-4}. Assume that $m\geq 16$. By Theorem \ref{main theorem}, there exists a $3$-SCHGDD of type $(n,4^{m/4})$ for any $n\equiv 0,2\ ({\rm mod}\ 6)$ and $m\equiv 4\ ({\rm mod}\ 12)$. We construct the $3$-SCHGDD on $I_n\times Z_m$ with the group set $\{\{i\}\times Z_m:i\in I_n\}$ and the hole set $\{I_n\times \{j,m/4+j,m/2+j,3m/4+j\}:0\leq j\leq m/4-1\}$. Denote the set of its base blocks by ${\cal B}^*$. For each $i\in I_n$, construct a strictly cyclic $3$-GDD of type $4^{m/4}$ on $\{i\}\times Z_m$ with the group set $\{\{i\}\times \{j,m/4+j,m/2+j,3m/4+j\}:0\leq j\leq m/4-1\}$, which exists by Lemma \ref{3-sCGDD}. Denote the set of its base blocks by ${\cal C}^*_i$. By Lemma \ref{n-4}, we can construct an optimal $2$-D $(n\times 4,3,1)$-OOC with $J(n\times 4,3,1)=n(2n-1)/3$ codewords $\cal D$ on the set $I_n\times \{0,m/4,m/2,3m/4\}$. Now we construct the required $2$-D $(n\times m,3,1)$-OOC on $I_n\times Z_m$. It is readily checked that ${\cal B}^*\cup(\bigcup_{i\in I_n}{\cal C}^*_i)\cup{\cal D}$ is the set of the required $n(nm-2)/6$ codewords. \qed

Combining the results of Lemma \ref{relation} and Theorem \ref{OOC-1}, we have the following corollary.

\begin{Corollary}\label{OOC-2} There exists an optimal $2$-D $(n\times m,3,1)$-OOC with $J(n\times m,3,1)$ codewords for any $n\equiv 0,4\ ({\rm mod}\ 12)$ and $m\equiv 2\ ({\rm mod}\ 6)$, or $n\equiv 0,8\ ({\rm mod}\ 24)$ and $m\equiv 1\ ({\rm mod}\ 3)$.
\end{Corollary}

Finally we make a remark. It is known in \cite{wsy} that an optimal $2$-D $(n\times m,3,1)$-OOC with $J(n\times m,3,1)$ codewords exists for any positive integer $n$ and odd integer $m$ except for $m=1$ and $n\equiv 5\ ({\rm mod}\ 6)$, where in the latter case, an optimal $2$-D $(n\times 1,3,1)$-OOC only has $J(n\times 1,3,1)-1$ codewords. For the case of $m$ even, it seems that the problem is still open. Here we only intend to present one example to show how to make optimal $2$-D OOCs via SCHGDDs. Actually more new optimal $2$-D $(n\times m,3,1)$-OOCs can be obtained by using $3$-SCHGDDs; especially, after some counting arguments, one can prove that in some cases, the upper bound of a $2$-D $(n\times m,3,1)$-OOC is Johnson bound minus one, and using similar techniques as above, some infinite families of optimal $2$-D OOCs attaining Johnson bound minus one can be obtained. This problem merits further study.

\subsection*{Acknowledgements}
The authors would like to thank the two anonymous referees for their
helpful comments and valuable suggestions. Research of this paper was partially carried out while the first author was visiting University of Victoria. He expresses his sincere thanks to China Scholarship Council for financial support and to the Department of Mathematics and Statistics in University of Victoria for their kind hospitality.


\begin{thebibliography}{Z}

\bibitem{aa}
R. J. R. Abel and A. Assaf, {\it Modified group divisible designs with block size $5$ and $\lambda=1$}, Discrete Math., 256(2002), 1-22.

\bibitem{abg}
R. J. R.  Abel, F. E. Bennett and M. Greig, PBD-closure, in: {\it
CRC Handbook of Combinatorial Designs} (C. J. Colbourn and J. H.
Dinitz, eds.), CRC Press, (2007), 247-255.

\bibitem{am1}
T. L. Alderson and K. E. Mellinger, {\it $2$-Dimensional optical
orthogonal codes from Singer groups}, Discrete Appl. Math.,
157(2009), 3008-3019.

\bibitem{a}
A. Assaf, {\it Modified group divisible designs}, Ars Combin., 29(1990), 13-20.

\bibitem{a2}
A. Assaf, {\it An application of modified group divisible designs}, J. Combin. Theory (A), 68(1994), 152-168.

%\bibitem{a1}
%A. Assaf, {\it Modified group divisible designs with block size $4$ and $\lambda>1$}, Australas. J. Combin., 16(1997), 229-238.

\bibitem{aw}
A. Assaf and R. Wei, {\it Modified group divisible designs with block size $4$ and $\lambda=1$}, Discrete Math., 195(1999), 15-25.

\bibitem{bjl}
T. Beth, D. Jungnickel and H. Lenz, {\it Design Theory}, Cambridge
University Press, Cambridge, UK, 1999.

\bibitem{bcrw}
D. Bryant, Y. Chang, C. A. Rodger, and R. Wei, {\it Two-dimensional balanced sampling plans excluding contiguous units}, Commun. Statist. Theory Methods, 31(2002), 1441-1455.

%\bibitem{bsh}
%A. E. Brouwer, A. Schrijver, and H. Hanani, {\it Group divisible designs
%with block size four}, Discrete Math., 20(1977), 1-10.

\bibitem{cww}
H. Cao, L. Wang, and R. Wei, {\it The existence of HGDDs with block size four and its application to double frames}, Discrete Math., 309(2009), 945-949.


\bibitem{cgz}
K. Chen, G. Ge, and L. Zhu, {\it Starters and related codes}, J. Statist. Plann. Inference, 86(2000), 379-395.

\bibitem{czl}
Z. Chen, {\it The existence of balanced difference families and perfect difference families}, Master degree thesis, Guangxi Normal University, 2008.

\bibitem{csw}
F. R. K. Chung, J. A. Salehi, and V. K. Wei, {\it Optical orthogonal
codes: design, analysis and applications}, IEEE Trans. Inform.
Theory, 35(1989) 595-604.

%\bibitem{cd}
%C. J. Colbourn and J. H. Dinitz, {\it Generating sets in Steiner triple systems}, Mathematica Slovaca, 50(2000), 259-269.

\bibitem{cl}
C. J. Colbourn and A. C. H. Ling, {\it A class of partial triple system with applications in survey sampling}, Commun. Statist. Theory Methods, 27(1998), 1009-1018.

\bibitem{dx}
C. Ding and C. Xing, {\it Several classes of $(2^m-1,w,2)$ optical orthogonal codes}, Discrete Appl. Math., 128(2003), 103-120.

\bibitem{d}
J. H. Dinitz, {\it Starters}, in: {\it CRC Handbook of Combinatorial Designs} (C. J. Colbourn and J. H. Dinitz, eds.), CRC Press, (2007), 622-628.

\bibitem{fwc}
T. Feng, X. Wang, and Y. Chang, {\it Semi-cyclic holey group divisible designs with block size three}, Designs, Codes and Cryptography, DOI 10.1007/s10623-013-9859-7.

\bibitem{fm}
R. Fuji-Hara and Y. Miao, {\it Optimal orthogonal codes: Their
bounds and new optimal constructions}, IEEE Trans. Inform. Theory,
46(2000), 2396-2406.

\bibitem{jiang}
R. P. Gallant, Z. Jiang, and A. C. H. Ling, {\it The spectrum of cyclic
group divisible designs with block size three}, J. Combin. Designs, 7(1999), 95-105.

%\bibitem{gr}
%G. Ge and R. S. Rees, {\it On group-divisible designs with group size four and group-type $g^um^1$}, Designs, Codes and Cryptography, 27(2002), 5-24.

%\bibitem{grz}
%G. Ge, R.S. Rees, and L. Zhu, {\it Group-divisible designs with block size four and group-type $g^um^1$ with $m$ as large as or as small as possible}, J. Combin. Theory (A), 98(2002), 357-376.

\bibitem{gww}
G. Ge, J. Wang, and R. Wei, {\it MGDD with block size $4$ and its application to sampling designs}, Discrete Math., 272(2003), 277-283.

\bibitem{gw}
G. Ge and R. Wei, HGDDs with block size four, Discrete Math., 279(2004), 267-276.

\bibitem{hrs}
A. S. Hedayat, C. R. Rao, and J. Stufken, {\it Sampling plans excluding contiguous units}, J. Statist. Plann. Inference, 19(1988), 159-170.

\bibitem{lc}
A. C. H. Ling and C. J. Colbourn, {\it Modified group divisible designs with block size four}, Discrete Math., 219(2000), 207-221.

%\bibitem{mr}
%R. Mathon and A. Rosa, $2$-$(v,k,\lambda)$ designs of small orders, in: {\it
%CRC Handbook of Combinatorial Designs} (C. J. Colbourn and J. H.
%Dinitz, eds.), CRC Press, (2006), 25-58.

\bibitem{sb}
J. A. Salehi and C. A. Brackett, {\it Code division multiple-access
techniques in optical fiber networks-Part II: Systems performance
analysis}, IEEE Trans. Communications, 37(1989), 834-842.

\bibitem{s}
D. R. Stinson, {\it A general construction for group divisible designs}, Discrete Math., 33(1981), 89-94.

\bibitem{sz}
D. R. Stinson and L. Zhu, {\it On the existence of MOLS with equal-sized holes}, Aequationes Math., 33(1987), 96-105.

\bibitem{stu}
J. Stufken, {\it Combinatorial and statistical aspects of sampling plans to avoid the selection of adjacent units}, J. Combin. Inform. System Sci., 18(1993), 81-92.

\bibitem{wsy}
J. Wang, X. Shan, and J. Yin, {\it On constructions for optimal
two-dimentional optical orthogonal codes}, Designs, Codes and
Cryptography, 54(2010), 43-60.

\bibitem{wy}
J. Wang and J. Yin, {\em Two-dimensional optical orthogonal codes
and semicyclic group divisible designs}, IEEE Trans. Inform. Theory,
56(2010), 2177-2187.

\bibitem{wlc}
L. Wang and Y. Chang, {\em Combinatorial constructions of optimal $3$-D $(u\times v\times w,3,1)$-OOCs}, submitted.


\bibitem{wc}
X. Wang and Y. Chang, {\it The spectrum of $(gv,g,3,\lambda)$-DF in $Z_{gv}$}, Science in China (A), 52(2009), 1004-1016.

\bibitem{wcf}
X. Wang, Y. Chang, and T. Feng, {\it Optimal $2$-D $(n\times m,3,2,1)$-optical orthogonal codes}, IEEE Trans. Inform. Theory, 59(2013), 710-725.

\bibitem{wei}
R. Wei, {\it Group divisible designs wih equal-sized holes}, Ars Combin., 35(1993), 315-323.

\bibitem{wei1}
R. Wei, {\it Cyclic BSEC of block size $3$}, Discrete Math., 250(2002), 291-298.

\bibitem{wr}
J. H. Wright, {\it Two-dimensional balanced sampling plans excluding adjacent units}, J. Statist. Plann. Inference, 138(2008), 145-153.

\bibitem{yk}
G. C. Yang and W. C. Kwong, {\it Performance comparison of
multiwavelength CDMA and WDMA+CDMA for fiber-optic networks,} IEEE
Trans. Communications, 45(1997), 1426-1434.

\bibitem{zhu}
L. Zhu, {\it Some recent developments on BIBDs and related designs}, Discrete Math., 123(1993), 189-214.

\end{thebibliography}
\end{document}